\documentclass[11pt]{article}

\title{Adaptive Thermostats for Noisy Gradient Systems}
\author{Benedict Leimkuhler and Xiaocheng Shang\footnote{Corresponding author. Email: \href{mailto:x.shang@ed.ac.uk}{x.shang@ed.ac.uk}} \\
\small{School of Mathematics, University of Edinburgh, Edinburgh, EH9 3FD, UK} }

\date{\today}

\usepackage{cite}

\usepackage{hyperref}
\usepackage{graphicx}

\usepackage{epstopdf} 

\usepackage{caption}
\DeclareGraphicsExtensions{.eps,.mps,.pdf,.jpg,.png}
\graphicspath{{figures/}{../figures/}}

\usepackage{fancyhdr}

\usepackage{amsmath}
\usepackage{bm}

\usepackage{amsfonts}

\usepackage{dsfont}

\usepackage[titletoc,toc,title]{appendix}

\usepackage[margin=3.0cm]{geometry}

\newcommand{\Tr}{\mathrm{Tr}}
\newcommand{\E}{\mathbb{E}}
\newcommand{\Var}{\mathrm{Var}}
\newcommand{\Cov}{\mathrm{Cov}}
\usepackage[usenames]{color} 
\usepackage{amssymb} 
\usepackage{amsmath} 
\usepackage{amsthm} 
\usepackage[utf8]{inputenc} 
\usepackage{bm}
\usepackage{mathrsfs}

\def\dd{\mathrm{d}}

\def\curlyL{\mathscr{L}}

\newcommand{\be}{\begin{equation}}
\newcommand{\ee}{\end{equation}}
\newcommand{\bea}{\begin{eqnarray}}
\newcommand{\eea}{\end{eqnarray}}

\renewcommand{\vec}[1]{ {\bm #1}}

\newcommand{\y}{{\vec{y}}}
\newcommand{\q}{{\vec{q}}}
\newcommand{\p}{{\vec{p}}}
\newcommand{\F}{{\vec{F}}}
\newcommand{\M}{{\vec{M}}}
\newcommand{\I}{{\vec{I}}}
\newcommand{\thetaB}{{\vec{\theta}}}

\usepackage{mathrsfs}
\def\curlyL{\mathscr{L}}
\usepackage{array}
\newcolumntype{C}[1]{>{\centering\let\newline\\\arraybackslash\hspace{0pt}}m{#1}}

\usepackage{soul}
\usepackage{color}

\begin{document}

\maketitle

\begin{abstract}
  We study numerical methods for sampling probability measures in high dimension where the underlying model is only approximately identified with a gradient system. Extended stochastic dynamical methods are discussed which have application to multiscale models, nonequilibrium molecular dynamics, and Bayesian sampling techniques arising in emerging machine learning applications. In addition to providing a more comprehensive discussion of the foundations of these methods, we propose a new numerical method for the adaptive Langevin/stochastic gradient Nos\'{e}--Hoover thermostat that achieves a dramatic improvement in numerical efficiency over the most popular stochastic \mbox{gradient} methods reported in the literature. We also demonstrate that the newly established method inherits a superconvergence property (fourth order convergence to the invariant measure for configurational quantities) recently demonstrated in the setting of Langevin dynamics. Our findings are verified by numerical experiments.
\end{abstract}

\pagenumbering{arabic}
\section{Introduction}

Stochastic thermostats~\cite{Samoletov2007,Leimkuhler2009,Samoletov2011} are powerful tools for sampling probability measures on high-dimensional spaces.   These methods combine an extended dynamics with degenerate stochastic perturbation to ensure ergodicity.    The traditional use of thermostats in molecular dynamics is to sample a well-specified equilibrium system involving a known force field which is the gradient of a potential energy function.   Recently, however, these techniques have become increasingly popular for problems of more general form, including the following:
\begin{itemize}
  \item multiscale models in which the forces are obtained by approximate sampling in another scale regime~\cite{Praprotnik2005,Praprotnik2008,Fedosov2009,Hijon2010,Li2014a,Mones2014};
  \item nonequilibrium physical models in which the potential energy function either is evolving or does not completely specify the system~\cite{Soddemann2003,Keaveny2005,Pastorino2007,Lisal2007,Peters2012,Stella2014};
  \item Bayesian machine learning applications in which a dataset defines an objective function which leads to an effective force law~\cite{Welling2011,Ahn2012,Patterson2013,Chen2014,Ding2014,Vollmer2015,Shang2015}.
\end{itemize}
In this article, we consider thermostats and numerical methods for sampling an underlying probability measure in the presence of error, under the assumption that the errors are random with a simple distributional form and unknown, but constant or slowly varying,  parameters.  In the cases considered, these methods are simple to implement, robust, and efficient.

\subsection{Thermostats}

The main tool that we employ in this article is the general concept of a thermostat as a (stochastic) distributional control for a dynamical system.  These methods originate in molecular dynamics, and it is simplest to explain them in that context.   Classical molecular dynamics tracks the motion of individual atoms determined by Newton's law in the microcanonical ($NVE$) ensemble, where energy (i.e., the Hamiltonian of the system) is always conserved~\cite{Allen1989,Frenkel2001,Leimkuhler2005,Hairer2006}. However, constant energy is not the appropriate setting of a real-world laboratory environment. In most cases, one wishes instead to sample the canonical ($NVT$) ensemble, where temperature, as an intensive variable, is conserved, by using thermostat techniques~\cite{Frenkel2001,Hunenberger2005}.

The idea of a thermostat is to modify dynamics so that a prescribed invariant measure is sampled.  There are competing aims in this type of work.  For example, one may wish to perturb the underlying Newtonian dynamics minimally, so that temporal correlations are preserved, or one may be interested in sampling rare events in a system with metastable states; thus a variety of methods have been developed.    The most obvious proposals, and also the oldest, are Brownian and Langevin dynamics.   In Brownian (sometimes called ``overdamped Langevin'') dynamics, the system is
\begin{equation} \label{eq:Brownian}
  {\rm d}\q = -\lambda \nabla U(\q){\rm d} t + \sqrt{2\beta^{-1} \lambda} {\rm d {\bf W}} \, ,
\end{equation}
where $\q$ represents a $3N$-dimensional vector of time-dependent random variables, ${\rm d {\bf W}}$ represents a vector of infinitesimal Wiener increments, $\beta$ is a positive parameter (proportional to the reciprocal temperature), $U$ is the potential energy function, and $\lambda$ is a free parameter which represents a time-rescaling.  It can be shown~\cite{Cances2007} that this system~\eqref{eq:Brownian} ergodically samples the Gibbs--Boltzmann probability distribution $\bar{\rho}_{\beta}  \propto  {\rm exp}(-\beta U)$.   For simplicity, we assume that the configurations $\q$ are restricted to a compact and simply connected domain $\Omega_{\q}$.   In molecular dynamics applications, the starting point is the potential energy function, which is usually assumed to be a semiempirical formula constructed from primitive functions via an a priori parameter fitting procedure.  Alternatively, one may assume that it is the probability distribution that is specified and that the potential energy is constructed from it via
\[
  U = -\beta^{-1} \ln \rho \, ,
\]
which, of course, requires that $\rho>0$.   In many applications it is found that the use of a first order dynamics such as~\eqref{eq:Brownian} is inefficient or introduces unphysical dynamical properties, and one employs, instead, the Langevin dynamics method:
\begin{align}
  {\rm d} \q & =  \M^{-1} {\p} {\rm d} t \, , \label{eq:Langevin-1} \\
  {\rm d} \p & =  -\nabla U(\q) {\rm d} t - \gamma {\p} {\rm d} t + \sqrt{2 \beta^{-1}\gamma}\M^{1/2} {\rm d} {\bf W} \, . \label{eq:Langevin-2}
\end{align}
Again, $\gamma$ in these equations is a free parameter, termed the ``friction constant".  It is related to the timescale on which the variables of the system interact with particles of a fictitious extended ``bath'', but it cannot be associated with a simple time-rescaling of the equations of motion and is thus different from $\lambda$ in~\eqref{eq:Brownian}.  It is a little more involved to show that~\eqref{eq:Langevin-1}--\eqref{eq:Langevin-2} ergodically~\cite{Mattingly2002} samples the distribution with density
$\rho_{\beta} \propto  {\rm exp}(-\beta H(\q, \p))$, where $H(\q, \p) = {\p}^T\M^{-1}{\p}/2 + U(\q)$.  In molecular dynamics, the $3N \times 3N$ matrix $\M$ is typically diagonal and contains the masses of atoms, ${\p}$ represents the momentum vector, and $H$ is the Hamiltonian or energy function.  In more general settings, the masses and friction coefficient may be treated as free parameters, and by computing long trajectories of~\eqref{eq:Langevin-1}--\eqref{eq:Langevin-2}, one may obtain averages with respect to $\bar{\rho}_{\beta}(\q)$; i.e., if $\{\left(\q(\tau),{\p}(\tau)\right): \tau \geq 0\}$ is a path generated by solving the SDE system~\eqref{eq:Langevin-1}--\eqref{eq:Langevin-2}, one has, for suitable test functions $\phi(\q)$, and under certain conditions on the potential energy function $U$~\cite{Mattingly2002},
\[
  \lim_{\tau\rightarrow \infty} \tau^{-1} \int_{0}^{\tau} \phi(\q(\tau)) \, {\rm d} \tau = \int_{\Omega_{\q}} \phi(\q) \bar{\rho}_{\beta}(\q) \, {\rm d } \omega_\q \, ,
\]
where ${\rm d}\omega_\q = {\rm d} \q_1 {\rm d} \q_2 \dots {\rm d} \q_N$.  In other words, the projected path defines a sampler for the density $\bar{\rho}_{\beta}$.

Langevin dynamics can thus be seen as an extended system which allows sampling to be performed in a reduced cross section of phase space by marginalization over long trajectories; this is the essential property of a thermostat.   Other types of thermostats include Nos\'{e}--Hoover--Langevin (NHL) dynamics~\cite{Samoletov2007,Leimkuhler2009} and various generalized schemes (see, e.g.,~\cite{Leimkuhler2010}).      In these methods, one adds additional auxiliary variables which are meant to control the dynamics (via a negative feedback loop), and the auxiliary variables are then further coupled to stochastic processes of Ornstein--Uhlenbeck type which can provide ergodicity~\cite{Leimkuhler2009}. (Note that the use of purely deterministic approaches, such as Nos\'{e}--Hoover, results in ergodicity issues~\cite{Legoll2006,Legoll2009}.)   The use of auxiliary variables can provide a degree of flexibility in the design of the thermostat, for example, allowing the treatment of systems arising in fluid dynamics~\cite{Dubinkina2010} or imposing an isokinetic constraint~\cite{Leimkuhler2013b}.  Very recently, we have further generalized the NHL method to obtain pairwise Nos\'{e}--Hoover--Langevin (PNHL), which is a momentum-conserving thermostat and thus applicable to the simulation of hydrodynamic behavior in complex fluids and polymers in mesoscales~\cite{Leimkuhler2015}.

\subsection{Noisy Gradients}

The gradient (or Hamiltonian) structure is essential to the nature of all the methods described above since it is only by use of this feature that the underlying Fokker--Planck equation can be shown to have the desired steady state solution.  However, in many applications, in particular multiscale modelling, the force is corrupted by significant approximation error and cannot be viewed as the gradient of a single global potential function.   One imagines a large extended system involving configurational variables $\q$ and $\y$, with $(\q, \y) \in \Omega_{\q}\times \Omega_{\y}$ (compact), and an overall distribution described by a Gibbs--Boltzmann density
\[
  \tilde{\rho}(\q, \y)= Z^{-1} {\rm exp} \left(-\beta \tilde{U} (\q,\y)\right) \, ,
\]
where $Z$ is a normalizing constant so that $\tilde{\rho}$ is a probability density. One calculates the mean force acting on $\q$, $\hat{f}(\q)$, by averaging the forces in the extended Gibbsian system, $\tilde{f}(\q,\y)$, as
\[
  \hat{f}(\q) = \int_{\Omega_\y} \tilde{f}(\q,\y) \tilde{\rho}(\q, \y) \, {\rm d}\omega_\y \, .
\]
If, as would typically be assumed, $\tilde{f}(\q,\y) = - \nabla_\q \tilde{U} (\q,\y)$, i.e., the force in the extended system is conservative, then we may interpret $\hat{f}$ as a conservative force as well, specifically the gradient of the potential of mean force, which is given by
\[
  \hat{U}(\q) = -\beta^{-1}  \ln \int_{\Omega_\y} {\rm exp}\left(-\beta \tilde{U}(\q, \y)\right) \, {\rm d}\omega_\y \, .
\]

The challenge arises when this integral must be approximated.   For example, if this is done by Monte Carlo integration, for fixed $\q$, one generates samples $\y^1, \y^2, \dots, \y^k$ from the distribution with density $\tilde{\rho}(\q,\y)$ and thus approximates
the mean force by
\[
  \bar{f}^k (\q) = k^{-1} \sum_{i=1}^k \tilde{f}\left(\q,\y^i\right) \, .
\]
In practice most systems constructed in this way, for example, those arising in mixed quantum and classical molecular models~\cite{Bornemann1996}, will admit very substantial errors in the forces; that is,
\[
  \bar{f}^k (\q) = \hat{f}(\q) + \Delta^{k}(\q) \, .
\]
Depending on the method of computation, it may be reasonable to assume that the errors $\Delta^k$ are normally distributed with zero mean, which is justified by the central limit theorem~\cite{Ash2012}, but the variance of the errors is generally not known and will be dependent on the location $\q$ where they are computed; thus we would expect
\begin{equation}\label{eq:Error_Sigma}
  \Delta^k(\q) \sim {\cal N}\left(\vec{0},\boldsymbol{\Sigma}^k(\q)\right) \, ,
\end{equation}
where $\boldsymbol{\Sigma}^k(\q)$ is an unknown covariance matrix. It should be noted that the assumption of the errors being Gaussian distributed is also often adopted in Bayesian inverse problems~\cite{Dashti2013} and elsewhere.

The most straightforward approach to the problem is to first treat the estimation problem for $\boldsymbol{\Sigma}^k$ separately, by some means, and then to use this within a standard Brownian or Langevin dynamics algorithm. The difficulty is that this requires a high level of local accuracy in the calculations, which is likely to be burdensome and involve redundant computation.  What we would prefer to do is to resolve the correct target distribution by a global calculation.

This problem has recently been encountered in the data science community, where it has attracted considerable attention~\cite{Welling2011,Ahn2012,Patterson2013,Chen2014,Ding2014,Vollmer2015,Shang2015}. To illustrate, we consider the problem of Bayesian sampling~\cite{Robert2004,Brooks2011}, where one is interested in correctly drawing states from a posterior probability density defined as
\begin{equation}\label{eq:posterior_Bayesian}
  \pi(\thetaB|\mathbf{X}) \propto \pi(\mathbf{X}|\thetaB) \pi(\thetaB) \, ,
\end{equation}
where $\thetaB$ is the parameter vector of interest, $\mathbf{X}$ represents the entire dataset, and, $\pi(\mathbf{X}|\thetaB)$ and $\pi(\thetaB)$ represent the likelihood and prior distributions, respectively.  In these applications, the distribution parameters are interpreted as the configuration variables ($\thetaB\equiv \q$).  We introduce a potential energy $U(\thetaB)$ by defining $\pi(\thetaB|\mathbf{X}) \propto \exp(-\beta U(\thetaB))$; thus taking the logarithm of~\eqref{eq:posterior_Bayesian} gives
\begin{equation}
  U(\thetaB) = -\log \pi(\mathbf{X}|\thetaB)-\log \pi(\thetaB) \, .
\end{equation}
Assuming the data are independent and identically distributed (i.i.d.), the logarithm of the likelihood distribution can then be calculated as
\begin{equation}
  \log \pi(\mathbf{X}|\thetaB) = \sum_{i=1}^N \log \pi(\mathbf{x}_{i}|\thetaB) \, ,
\end{equation}
where $N$ is the size of the entire dataset.

However, in machine learning applications, one often finds that directly sampling with the entire large-scale dataset is computationally infeasible. For instance, standard Markov chain Monte Carlo (MCMC) methods~\cite{Metropolis1953} require the calculation of the acceptance probability and the creation of informed proposals based on the whole dataset, while the gradient is evaluated through the whole dataset in the hybrid Monte Carlo (HMC) method~\cite{Duane1987,Horowitz1991,Brooks2011}, again resulting in severe computational complexity.

In order to improve the efficiency of simulation, the so-called stochastic gradient Langevin dynamics (SGLD) was recently proposed~\cite{Welling2011} based on using a random (and much smaller, i.e., $\tilde{n} \ll N$) subset to approximate the likelihood of the dataset for given parameters,
\begin{equation}
  \log \pi(\mathbf{X}|\thetaB) \approx \frac{N}{\tilde{n}}\sum_{i=1}^{\tilde{n}}\log \pi(\mathbf{x}_{r_{i}}|\thetaB) \, ,
\end{equation}
where $\{\mathbf{x}_{r_{i}}\}^{\tilde{n}}_{i=1}$ represents a random subset of $\mathbf{X}$. Overall, the ``noisy'' potential energy now can be written as
\begin{equation}
  \tilde{U}(\thetaB) = -\frac{N}{\tilde{n}}\sum_{i=1}^{\tilde{n}}\log \pi(\mathbf{x}_{r_{i}}|\thetaB) - \log \pi(\thetaB) \, ,
\end{equation}
with ``noisy'' force $\tilde{\F}(\thetaB) = -\nabla \tilde{U}(\thetaB)$.

\subsection{Sampling Methods for Noisy Gradients}

The challenge is to identify a method to compute samples distributed according to the Gibbs distribution $\rho(\q) = Z^{-1}\exp(-\beta U(\q))$, where the only available information is a stochastically perturbed force
$\tilde{\F}(\q)$ defined in the previous section.

In the original SGLD method, samples are generated by Brownian
dynamics,
\begin{equation}
  \q_{n+1}  =  \q_n + \Delta t_n \tilde{\F}(\q_n)  + \sqrt{2\beta^{-1} \Delta t_n} \mathbf{R}_n \, ,
\end{equation}
where $\mathbf{R}_n$ is a vector of i.i.d.\ standard normal random variables. It should be emphasized that $\Delta t_n$ is a sequence of stepsizes decreasing to zero~\cite{Welling2011}. Although a central limit theorem associated with the decreasing stepsize sequence was established by Teh et al.~\cite{Teh2014}, a fixed stepsize is often preferred in practice, which is the choice in this article as in Vollmer et al.~\cite{Vollmer2015}, where a modified SGLD (mSGLD) is introduced:
\begin{equation}
  \q_{n+1}  =  \q_n + \Delta t \tilde{\F}(\q_n)  + \sqrt{2\beta^{-1} \Delta t} \left( \I - \frac{\Delta t}{4} \Cov\tilde{\F}(\q_n) \right) \mathbf{R}_n \, ,
\end{equation}
where
\begin{equation}
  \Cov\tilde{\F}_{ij} = \E \left[ \left( \tilde{\F}_{i} - \E(\tilde{\F}_{i})\right) \left( \tilde{\F}_{j} - \E(\tilde{\F}_{j}) \right)^{T} \right]
\end{equation}
is the covariance matrix of the noisy force.

A stochastic gradient Hamiltonian Monte Carlo (SGHMC) method was also proposed very recently by Chen et al.~\cite{Chen2014}, which incorporates a parameter-dependent diffusion matrix $\boldsymbol{\Sigma}(\q)$ (i.e., the covariance matrix of the noisy force).  $\boldsymbol{\Sigma}(\q)$ is intended to effectively offset the stochastic perturbation of the gradient. However, it is very difficult to accommodate $\boldsymbol{\Sigma}(\q)$ in practice; moreover, as pointed out in~\cite{Ding2014}, poor estimation of it may have a significant adverse influence in correctly sampling the target distribution unless the stepsize is small enough.

These problems challenge the conventional mechanism of thermostats.   An article of Jones and Leimkuhler~\cite{Jones2011} provides an alternative means of tackling this problem by showing that Nos\'{e}--Hoover dynamics is able to adaptively dissipate excess heat pumped into the system while maintaining the Gibbs (canonical) distribution.  In the setting of systems involving a driving stochastic perturbation, the adaptive Nos\'{e}--Hoover method is referred to as Ad-NH, with similar generalizations of Nos\'{e}--Hoover--Langevin (Ad-NHL) and Langevin dynamics (Ad-Langevin) available.    An idea equivalent to Ad-Langevin was very recently applied in the setting of Bayesian sampling for use in data science calculations by Ding et al.~\cite{Ding2014}, which they referred to as the stochastic gradient Nos\'{e}--Hoover thermostat (SGNHT).  It showed significant advantages over alternative techniques such as SGHMC~\cite{Chen2014}. However, the numerical method used by Ding et al.~\cite{Ding2014} is not optimal, neither in terms of its accuracy (measured per unit work) nor its stability (measured by the largest usable stepsize).

Although extended systems have been increasingly popular in molecular simulations, the mathematical analysis of the order of convergence, specifically in terms of the bias in averaged quantities computed using numerical trajectories, is not fully understood.  Using a splitting approach, we propose in this article an alternative numerical method for Ad-Langevin simulation that substantially improves on the existing schemes in the literature in terms of accuracy, robustness, and overall numerical efficiency.

The rest of the article is organized as follows. In Section~\ref{sec:Adaptive_Thermostats}, we describe the construction of adaptive formulations for noisy gradients including the Ad-Langevin/SGNHT method.  Section~\ref{sec:Numerical_Methods} considers the construction of numerical methods for solving the SDEs.   Numerical experiments are performed in Section~\ref{sec:Numerical_Experiments}. Our experiments are of a more limited nature in comparison with those of Ding et al.~\cite{Ding2014}, but we believe them to be representative of performance on a significant class of problems.   Finally, we summarize our findings in Section~\ref{sec:Conclusions}.

\section{Adaptive Thermostats for Noisy Gradients}
\label{sec:Adaptive_Thermostats}

In this section, we discuss the construction of thermostats to approximate samples with respect to the target measure (i.e., the correct marginalized Gibbs density) if the covariance matrix of the noisy force is constant, i.e., $\boldsymbol{\Sigma}(\q)=\sigma^{2}\I$ ($\sigma$ is a constant positive quantity). The procedure was outlined in the paper of Jones and Leimkuhler~\cite{Jones2011} and relies on the fact that a fixed amplitude noise perturbation engenders a shift of the auxiliary variable in the extended stationary distribution associated with the Nos\'{e}--Hoover thermostat.

If the system is not coming from a Newtonian dynamics model, then it is unclear that we need to rely on second order dynamics for this purpose.  To see why this is the case, we explain  what goes wrong if we try to use first order dynamics. In what follows, we assume that the covariance matrix of the noisy force is constant, although we ultimately intend to apply the method more generally (see recent work on a novel covariance-controlled adaptive Langevin thermostat that can handle parameter-dependent noise in~\cite{Shang2015}).  Even in the constant $\sigma$ case it is a nontrivial problem to extract statistics related to a particular target temperature, since we do not assume that $\sigma$ is known.

For $\sigma$ constant, let us first consider the SDE
\begin{align}
  \dd \q & = -\xi \nabla U(\q) \dd t + \sigma \dd {\bf W} \, , \label{eq:adBD1}\\
  \dd \xi & = \chi(\q) \dd t \label{eq:adBD2}
\end{align}
and seek $\chi(\cdot)$ so that an extended Gibbs distribution  with density of the form $\psi (\q,\xi)= \bar{\rho}_{\beta}(\q)\varphi(\xi)$ is (ergodically) preserved.  The variable $\xi$ is an auxiliary variable.  We do not generally care what its distribution is, but it is crucial that
\begin{enumerate}
  \item[(i)] the overall density is in product form, and
  \item[(ii)] $\varphi(\xi)\geq 0$ is normalizable and of a simple, easily sampled form.
\end{enumerate}
These conditions ensure that we can easily average out over the auxiliary variable to compute the averages of functions of $\q$ which are of greatest interest.\\

\par\noindent {\bf Proposition 1.} \emph{Let $\chi(\q) = -\beta^{-1}\Delta U(\q) +  \| \nabla U(\q)\|^2$; then~\eqref{eq:adBD1}--\eqref{eq:adBD2} preserves the modified Gibbs distribution
\[
  \tilde{\rho} (\q,\xi)= \bar{\rho}_{\beta}(\q) e^{-\beta (\xi-\hat{\gamma})^2/2} \, ,
\]
where $\hat{\gamma} =\beta \sigma^2/2$.}

\par\noindent {\emph{Proof.}}
The Fokker--Planck equation corresponding to~\eqref{eq:adBD1}--\eqref{eq:adBD2} is
\[
  \rho_t = {\cal L}^{\dagger} \rho: =  \xi \nabla\cdot \left(\nabla U (\q)\rho(\q,\xi)\right) + \frac{\sigma^2}{2} \Delta \rho - \frac{\partial}{\partial \xi}(\chi(\q) \rho) \, .
\]
Just insert $\tilde{\rho}$ into the operator ${\cal L}^{\dagger}$ to see that it vanishes.
\hfill $\Box$\\

Proposition 1 tells us that if we can solve system~\eqref{eq:adBD1}--\eqref{eq:adBD2}, under an assumption of ergodicity, we can compute averages with respect to the target Gibbs distribution without actually knowing the value of $\sigma$. $\sigma$ could be observed retrospectively by simply averaging $\xi$ during simulation, since $\langle \xi \rangle = \beta \sigma^2/2$.

The problem is that the dynamics~\eqref{eq:adBD1}--\eqref{eq:adBD2} is not quite what we want.
A typical numerical method for this system might be constructed based on modification of the Euler--Maruyama method:
\begin{align}
  \q_{n+1} & = \q_n -\Delta t \xi_n \nabla U(\q_n)+ \sigma \sqrt{\Delta t} \mathbf{R}_n \, ,  \label{eq:em1}\\
  \xi_{n+1} & = \xi_n + \Delta t \chi(\q_n) \, ; \label{eq:em2}
\end{align}
however, observe that this method requires separate knowledge of $\nabla U(\q)$ and $\sigma$, which is generally impossible a priori, as we assume that the force is polluted by unknown noise.   The form of the equations means that we evaluate the product of $\xi$ and the deterministic force, on the one hand, and the random perturbation, on the other hand, separately, and these contributions are independently scaled by $\Delta t$ and $\sqrt{\Delta t}$, respectively.

\subsection{The Adaptive Langevin (Ad-Langevin) Thermostat }

To adaptively control the invariant distribution, we consider the following second order formulation, which was first introduced in the paper of Jones and Leimkuhler~\cite{Jones2011}:
\begin{equation}
  \label{eq:Ad-L}
  \begin{aligned}
    \dd \q &= \M^{-1}\p\dd t \, , \\
    \dd \p &= \tilde{\F}(\q)\dd t - \xi \p\dd t + \sigma_{\mathrm{A}} \M^{1/2} \dd \mathbf{W}_{\rm A} \, , \\
    \dd \xi        &= {\mu}^{-1} \left[ \p^{T}\M^{-1}\p - N_{\mathrm{d}}k_{\mathrm{B}}T \right]\dd t \, .
  \end{aligned}
\end{equation}
In these equations, $\tilde{\F}(\q)$ is meant to represent a noisy gradient which may be thought of as being defined by the relation
\begin{equation}
\tilde{\F} (\q)  = -\nabla U(\q)  + \sigma \M^{1/2} \mathbf{R} \, ,
\end{equation}
where $\mathbf{R}=\mathbf{R}(t)$ is a collection of independent Gaussian white noise processes, i.e., $\langle\mathrm{R}_{i}(t)\mathrm{R}_{j}(s)\rangle = \delta_{ij} \delta(t-s)$, where $\delta_{ij}$ is the Kronecker delta and $\delta(t-s)$ is the Dirac delta function. $\sigma_{\rm A} \M^{1/2} \dd\mathbf{W}_{\rm A}$ indicates the artificial noise added into the system to enhance the ergodicity; i.e., the constant $\sigma_{\rm A}$ is known a priori.  All the components of the Wiener process ${\mathbf{W}}_{\rm A}(t)$ are assumed to be independent. $N_{\mathrm{d}}$ denotes the number of degrees of freedom of the system. $\mu$ is a coupling parameter which is referred to as the ``thermal mass''. $k_{\mathrm{B}}$ and $T$, satisfying the relation $\beta^{-1}=k_{\mathrm{B}}T$, represent the Boltzmann constant and system temperature, respectively.

A similar system (SGNHT) was used by Ding et al.~\cite{Ding2014}, who also explored its application to three examples from machine learning.  These experiments demonstrated that Ad-Langevin has superior performance compared to SGHMC in various applications, confirming the importance of adaptively dissipating additional noise in sampling. However, there remain two important  issues that we wish to address in this article: (1) the underlying dynamics of the Ad-Langevin method is not clear due to the presence of the stochastically perturbed gradient; (2) little attention has been paid to the design of optimal numerical methods for implementing Ad-Langevin with attention to stability and numerical efficiency.

One may wonder why the artificial noise is needed (i.e., $\sigma_{\rm A} \neq 0$), since we are assuming the presence of noise in the gradient itself.    The reason is as follows: in defining a numerical method for the noisy gradient system, the force (including the random perturbation) will in general be multiplied by $\Delta t$, where $\Delta t$ is the timestep.   On the other hand, the It\={o} rule implies that the scaling of random perturbations in an SDE should be by a factor proportional to $\sqrt{\Delta t}$; thus, effectively, if we are to relate the thermostatted method to a standard SDE, the standard deviation of the noise is reduced by multiplication by the factor $\sqrt{\Delta t}$.  The noise perturbation introduced at each timestep (and the effective diffusion) is thus reduced for small stepsizes and it is therefore important to inject additional artificial noise in order to stabilize the invariant distribution.    A rewriting of the Ad-Langevin system as a standard It\={o} SDE system makes clear the relation between the different terms
\begin{equation}
  \label{eq:Ad-L-2}
  \begin{aligned}
    \dd \q &= \M^{-1}\p\dd t \, , \\
    \dd \p &= -\nabla U(\q)\dd t + \sigma \sqrt{\Delta t} \M^{1/2} \dd {\mathrm{{\bf W}}} - \xi \p\dd t + \sigma_{\mathrm{A}} \M^{1/2} \dd\mathrm{{\bf W}}_{\rm A} \, , \\
    \dd \xi        &= {\mu}^{-1} \left[ \p^{T}\M^{-1}\p - N_{\mathrm{d}}k_{\mathrm{B}}T \right]\dd t \, ,
  \end{aligned}
\end{equation}
where ${\mathrm{{\bf W}}}={\mathrm{{\bf W}}}(t)$ is an additional vector of standard Wiener processes.

Let us note the main features of the dynamics~\eqref{eq:Ad-L-2}:
\begin{enumerate}
  \item[(i)] The equations are a combination of Langevin dynamics and Nos\'{e}--Hoover dynamics. If $\xi$ is constant in the equation for the momentum, then the system reduces to Langevin dynamics.  In the absence of noise, $\sigma_\mathrm{A}=0$ (and $\sigma=0$); then the system reduces to Nos\'{e}--Hoover.  The system~\eqref{eq:Ad-L-2} may be regarded as a sort of Langevin dynamics where the friction coefficient, rather than being fixed a priori, is automatically and adaptively determined in order to achieve the desired temperature (which is specified in the control law defining the evolution of $\xi$).
  \item[(ii)] The invariant distribution for the given system may be directly obtained by study of its Fokker--Planck equation.  Following~\cite{Jones2011}, it is straightforward to show that~\eqref{eq:Ad-L-2} has the following invariant distribution:
      \begin{equation}\label{eq:Invariant_Dist_Modified}
        \tilde{\rho}_{\beta}(\q,\p,\xi) = \frac{1}{Z}\exp\left({-\beta H(\q,\p)}\right)\exp\left( -\frac{\beta\mu}{2} (\xi - \hat{\gamma})^{2} \right) \, ,
      \end{equation}
      where $Z$ is the normalizing constant and
      \begin{equation}\label{eq:xi_bar}
        \hat{\gamma} = \frac{ \beta \left( \sigma^2_{\mathrm{F}} + \sigma^2_{\mathrm{A}} \right)}{2}  \, ,
      \end{equation}
      where $\sigma_{\rm F} = \sigma\sqrt{\Delta t}$.
      Observe that this means that if $\sigma_{\rm A}=0$, then, as ${\rm lim}_{\Delta t\rightarrow 0} \sigma_{\rm F}=0$, we find that $\xi$ tends to a variable which is normally distributed with mean zero.    Alternatively, if  $\sigma_{\rm A}\neq 0$, one would obtain
      \[
      \xi   \overset{\curlyL}{\to} {\cal N} \left( \frac{ \beta \sigma^2_{\mathrm{A}}}{2} , \beta^{-1}\mu^{-1} \right), \hspace{0.2in}  t \rightarrow \infty \, , \hspace{0.2in}
      \Delta t \to 0 \, ,
      \]
      where $\beta^{-1}\mu^{-1}$ is the variance and the symbol $ \overset{\curlyL}{\to} $ indicates that $\xi$ converges in probability law to a normally distributed random variable with the indicated parameters.  The order of the limits here is important: $t\rightarrow \infty$ first (to reach the invariant distribution), then $\Delta t\rightarrow 0$.
  \item[(iii)] The ergodicity of \eqref{eq:Ad-L-2} with respect to the distribution indicated above can easily be demonstrated by reference to H\"ormander's condition for hypoellipticity following the method in \cite{Mattingly2002}, as for Langevin dynamics.   The only additional step is to verify that the noise propagates into the $\xi$ variable, which follows due to its strong coupling to the momenta.
  \item[(iv)] This dynamics is a bit unusual in that it must be viewed as stepsize dependent, although we mention that such mixed systems are used in the study of backward error analysis~\cite{Leimkuhler2005}.   One simply thinks of the characteristics of stochastic paths associated with~\eqref{eq:Ad-L-2} as being stepsize dependent.  Although \eqref{eq:Ad-L-2} takes on the appearance of a standard It\={o} SDE system, we must bear in mind that in discretizing these equations the conservative force $\F(\q)$ and the associated noise term $\sigma \sqrt{\Delta t} \M^{1/2} \dd {\mathrm{{\bf W}}}$ must be evaluated together at every stage, since the formulation (\ref{eq:Ad-L-2}) is a notational device to make clear the properties of the system.
\end{enumerate}

\section{Numerical Methods for Adaptive Thermostats}
\label{sec:Numerical_Methods}

Since stochastic systems in most of the cases cannot be solved ``exactly'', splitting methods are often adopted in practice. For instance here, the vector field of the Ad-Langevin/SGNHT~\eqref{eq:Ad-L} can be split into four pieces which are denoted as ``A'', ``B'', ``O'', and ``D'', in such a way that each piece can be solved ``exactly'',
\begin{equation*}\label{eq:Splitting_SGNHT}
  \dd \left[ \begin{array}{c} \q \\ \p \\ \xi \end{array} \right] = \underbrace{\left[ \begin{array}{c} {\M}^{-1}\p \\ \vec{0} \\ 0 \end{array} \right] \dd t}_\mathrm{A} + \underbrace{\left[ \begin{array}{c} \vec{0} \\  -\nabla U(\q)  + \sigma \M^{1/2} \mathbf{R} \\ 0 \end{array} \right] \dd t }_\mathrm{B} + \underbrace{\left[ \begin{array}{c} \vec{0} \\ - \xi \p\dd t + \sigma_{\mathrm{A}} \M^{1/2} \dd\mathrm{{\bf W}}_{\rm A} \\ 0 \end{array} \right] }_\mathrm{O} + \underbrace{\left[ \begin{array}{c} \vec{0} \\ \vec{0} \\ G(\p) \end{array} \right] \dd t}_\mathrm{D} \, ,
\end{equation*}
where $G(\p) = {\mu}^{-1} \left[ \p^{T}\M^{-1}\p - N_{\mathrm{d}}k_{\mathrm{B}}T \right]$.

Clearly parts ``A'' and ``D'' can be solved ``exactly''.   As mentioned previously, the underlying dynamics for ``B'' is
\begin{equation}\label{eq:SDE_Piece_B}
  \dd \p = -\nabla U(\q)\dd t + \sigma_{\mathrm{F}} \M^{1/2} \dd {\mathrm{{\bf W}}} \, ,
\end{equation}
where $\q$ is fixed and $\sigma_{\mathrm{F}} = \sigma \sqrt{\Delta t}$.
Integrating~\eqref{eq:SDE_Piece_B} from 0 to $\Delta t$ gives the exact solution in distribution of this part as
\begin{align*}
  \p(\Delta t) & = \p(0) - \Delta t \nabla U(\q) +  \sqrt{\Delta t} \sigma_{\mathrm{F}} \M^{1/2} {\mathrm{{\bf R}}} \\
  & = \p(0) + \Delta t [-\nabla U(\q)  + \sigma \M^{1/2} \mathbf{R}] = \p(0) + \Delta t \tilde{\F}(\q) \, ,
\end{align*}
where ${\mathrm{{\bf R}}}$ is a vector of i.i.d.\ standard normal random variables. It should be noted that applying the Euler--Maruyama method to~\eqref{eq:SDE_Piece_B} gives the same result; thus, for constant force, Euler--Maruyama is ``exact''.

The ``O'' or ``Ornstein--Uhlenbeck'' part is usually stated with $\xi$ a positive constant, in which case the solution is found to be~\cite{Kloeden1992}
\begin{equation}\label{eq:Standard_OU_sol}
  \p(\Delta t) = e^{-\xi \Delta t} \p(0) + \sigma_{\mathrm{A}}\sqrt{ \frac{1-e^{-2\xi \Delta t}}{2\xi} }\M^{1/2}\mathbf{R} \, ,
\end{equation}
where $\mathbf{p}(0)$ is the initial value of the variable and $\mathbf{R}$ is a vector of i.i.d.\ standard normal random variables. However, the same formula (\ref{eq:Standard_OU_sol}) is easily seen to be valid for $\xi<0$, since the quantity $(1-e^{-2\xi \Delta t})/(2\xi)$ is strictly greater than zero unless $\xi=0$.  (The proof is obtained by following the standard procedure~\cite{Kloeden1992}.)  When $\xi=0$, one can simply replace $(1-e^{-2\xi \Delta t})/(2\xi)$ by its well-defined asymptotic limit,
\begin{equation}\label{eq:Standard_OU_sol_2}
  \p(\Delta t) = \p(0) + \sqrt{\Delta t} \sigma_{\mathrm{A}} \M^{1/2} \mathbf{R} \, .
\end{equation}

The generators associated with each piece are defined, respectively, as
\begin{equation*}\label{eq:Generators_Ad-Langevin}
  \begin{aligned}
    \mathcal{L}_\mathrm{A} &= \M^{-1}\p \cdot \nabla_{\q} \, , \\
    \mathcal{L}_\mathrm{B} &= -\nabla U(\q) \cdot \nabla_{\p} + \frac{{\sigma}^{2}_{\mathrm{F}}}{2}\Tr \left(\M\nabla^{2}_{\p}\right) \, , \\
    \mathcal{L}_\mathrm{O} &= -\xi \p \cdot \nabla_{\p} + \frac{{\sigma}^{2}_{\mathrm{A}}}{2} \Tr\left(\M\nabla^{2}_{\p}\right) \, , \\
    \mathcal{L}_\mathrm{D} &= G(\p) \frac{\partial}{\partial\xi} \, ,
  \end{aligned}
\end{equation*}
where $\sigma_{\mathrm{F}} = \sigma \sqrt{\Delta t}$ in part ``B'' is stepsize dependent.

Overall, the generator of the Ad-Langevin/SGNHT~\eqref{eq:Ad-L} system can be written as
\begin{equation}
  \mathcal{L} = \mathcal{L}_\mathrm{A} + \mathcal{L}_\mathrm{B} + \mathcal{L}_\mathrm{O} + \mathcal{L}_\mathrm{D} \, .
\end{equation}

The flow map (or phase space propagator) of the system can be written in the shorthand notation
\begin{equation*}
  { \cal F}_{t} = e^{t \mathcal{L} } \, ,
\end{equation*}
where the exponential map here denotes the solution operator.  Approximations of ${\cal F}_t$ can be obtained as products (taken in different arrangements) of exponentials of the splitting terms.   For example, the phase space propagation of the method proposed by Ding et al.~\cite{Ding2014} for the Ad-Langevin/SGNHT~\eqref{eq:Ad-L} system (denoted as ``SGNHT-N'') can be written as
\begin{equation}\label{eq:Splitting_SGNHT-N}
  \exp\left(\Delta t \hat{\mathcal{L}}_\mathrm{SGNHT-N} \right) = \exp\left(\Delta t \mathcal{L}_\mathrm{P}\right) \exp\left(\Delta t \mathcal{L}_\mathrm{A}\right) \exp\left(\Delta t \mathcal{L}_\mathrm{D}\right) \, ,
\end{equation}
where
\begin{equation}
  \mathcal{L}_\mathrm{P} = \mathcal{L}_\mathrm{B} + \mathcal{L}_\mathrm{O}
\end{equation}
and $\exp\left(\Delta t\mathcal{L}_f\right)$ represents the phase space propagator associated with the corresponding vector field $f$.  Because of its nonsymmetric structure, one anticipates first order convergence to the invariant measure (for any choice of $\sigma$).  Due to the naming of the component parts, the SGNHT-N method may be denoted by ``PAD''.

Overall, the SGNHT-N/PAD  integration method is as follows:
\begin{align*}
  \p_{n+1} &= \p_{n} + \Delta t \left( -\nabla U(\q_{n})  + \sigma \M^{1/2} \mathbf{R}'_{n} \right) - \Delta t \xi_{n} \p_{n} + \sqrt{\Delta t} \sigma_{\mathrm{A}} \M^{1/2} \mathrm{\bf R}_{n} \, , \\
  \q_{n+1} &=  \q_{n} + \Delta t \M^{-1}\p_{n+1} \, , \\
  \xi_{n+1} &= \xi_{n} +  \Delta t {\mu}^{-1} \left( \p^{T}_{n+1}\M^{-1}\p_{n+1} - N_{\mathrm{d}}k_{\mathrm{B}}T \right) \, ,
\end{align*}
where $\mathrm{\bf R}'_{n}$ and $\mathrm{\bf R}_{n}$ are vectors of i.i.d.\ standard normal random variables.

We propose symmetric alternative methods, such as the following
symmetric Ad-Langevin/\\SGNHT (SGNHT-S) splitting method:
\begin{equation}\label{eq:Splitting_SGNHT-S}
  e^{\Delta t\hat{\mathcal{L}}_\mathrm{SGNHT-S}} = e^{\frac{\Delta t}{2}\mathcal{L}_\mathrm{B}} e^{\frac{\Delta t}{2}\mathcal{L}_\mathrm{A}} e^{\frac{\Delta t}{2}\mathcal{L}_\mathrm{D}} e^{\Delta t\mathcal{L}_\mathrm{O}} e^{\frac{\Delta t}{2}\mathcal{L}_\mathrm{D}} e^{\frac{\Delta t}{2}\mathcal{L}_\mathrm{A}} e^{\frac{\Delta t}{2}\mathcal{L}_\mathrm{B}} \, ,
\end{equation}
where exact solvers for parts ``B'' and ``O'' derived above are applied.  The SGNHT-S method may be referred to as ``BADODAB'', where it should be noted that the various operations are symmetrically applied and the steplengths are uniform and span the interval $\Delta t$.   Other symmetric splittings are considered below.

The SGNHT-S numerical integration method may be written as
\begin{align*}
  \p_{n+1/3} &= \p_{n} + (\Delta t/2) \left( -\nabla U(\q_{n})  + \sigma \M^{1/2} \mathbf{R}'_{n} \right) \, , \\
  \q_{n+1/2} &=  \q_{n} + (\Delta t/2) \M^{-1}\p_{n+1/3} \, , \\
  \xi_{n+1/2} &= \xi_{n} +  (\Delta t/2) {\mu}^{-1} \left( \p^{T}_{n+1/3}\mathbf{M}^{-1}\p_{n+1/3} - N_{\mathrm{d}}k_{\mathrm{B}}T \right) \, , \\
  \mathrm{if}\: (\xi_{n+1/2} \neq 0):  \ \p_{n+2/3} &= e^{-\xi_{n+1/2} \Delta t}\p_{n+1/3} + \sigma_{\mathrm{A}} \sqrt{ (1-e^{-2\xi_{n+1/2} \Delta t})/(2\xi_{n+1/2}) } \M^{1/2} \mathrm{\bf R}_{n} \, , \\
  \mathrm{else}:  \ \p_{n+2/3} &= \p_{n+1/3} + \sqrt{\Delta t} \sigma_{\mathrm{A}} \M^{1/2} \mathbf{R}_{n} \, , \\
  \xi_{n+1} &= \xi_{n+1/2} +  (\Delta t/2) {\mu}^{-1} \left( \p^{T}_{n+2/3}\mathbf{M}^{-1}\p_{n+2/3} - N_{\mathrm{d}}k_{\mathrm{B}}T \right) \, , \\
  \mathbf{q}_{n+1} &=  \q_{n+1/2} + (\Delta t/2) \M^{-1}\p_{n+2/3} \, , \\
  \p_{n+1} &= \p_{n+2/3} + (\Delta t/2) \left( -\nabla U(\q_{n+1})  + \sigma \M^{1/2} \mathbf{R}'_{n+1} \right) \, .
\end{align*}
The force computed at the end of each timestep can be reused at the start of the next step; thus only one force calculation is needed in SGNHT-S at each timestep, the same as for SGNHT-N. In practice, one could replace the exponential and square root operations in the exact solver of the ``O'' part by their respective well-defined asymptotic expansions to reduce the computational cost.

\subsection{Order of Convergence of Ad-Langevin/SGNHT}
\label{subsec:Order_of_Convergence}

The analysis of the accuracy of ergodic averages (averages with respect to the invariant measure) in stochastic numerical methods can be performed using the framework of long-time Talay--Tubaro expansion, as developed in~\cite{Talay1990,Debussche2012,Leimkuhler2013,Leimkuhler2013a,Abdulle2014a,Abdulle2014,Leimkuhler2013c}.
In what follows we compare the order of convergence of the two Ad-Langevin/SGNHT methods with a clean gradient.

For a splitting method described by $ \mathcal{L}= \mathcal{L}_{\alpha} + \mathcal{L}_{\beta} + \dots + \mathcal{L}_{\zeta} $, we define the effective operator $\hat{\mathcal{L}}^{\dag}$ associated with the perturbed system obtained using the numerical method with stepsize $\Delta t$ by the relation
\[
\exp\left(\Delta t \hat{\mathcal{L}}^{\dag}\right) = \exp\left(\Delta t \mathcal{L}_{\alpha}^{\dag}\right) \exp\left(\Delta t\mathcal{L}_{\beta}^{\dag}\right) \dots \exp\left(\Delta t\mathcal{L}_{\zeta}^{\dag}\right) \, .
\]
This operator can be computed using the Baker--Campbell--Hausdorff (BCH) expansion and can thus be viewed as a perturbation of the exact Fokker--Planck operator $\mathcal{L}^{\dag}$:
\begin{equation}\label{eq:operator_perturbation}
  \hat{\mathcal{L}}^{\dag} = \mathcal{L}^{\dag} + \Delta t\mathcal{L}^{\dag}_{1} + \Delta t^{2}\mathcal{L}^{\dag}_{2} + O(\Delta t^{3})
\end{equation}
for some perturbation operators $\mathcal{L}^{\dag}_{i}$.

We also define the invariant distribution $\hat{\rho}$ associated with the numerical method as an approximation of the target invariant distribution $\tilde{\rho}_{\beta}$:
\begin{equation}\label{eq:distribution_perturbation}
  \hat{\rho} = \tilde{\rho}_{\beta}\left[ 1+\Delta tf_{1}+\Delta t^{2}f_{2}+\Delta t^{3}f_{3}+O(\Delta t^{4}) \right]
\end{equation}
for some correction functions $f_{i}$ satisfying $\langle f_{i} \rangle=0$.

Substituting $\hat{\mathcal{L}}^{\dag}$ and $\hat{\rho}$ into the stationary Fokker--Planck equation
\begin{equation*}
  \hat{\mathcal{L}}^{\dag}\hat{\rho} = 0
\end{equation*}
yields
\begin{equation*}
  \left( \mathcal{L}^{\dag} + \Delta t\mathcal{L}^{\dag}_{1} + \Delta t^{2}\mathcal{L}^{\dag}_{2} + O(\Delta t^{3}) \right)\left(\tilde{\rho}_{\beta}\left[1+\Delta tf_{1}+\Delta t^{2}f_{2}+\Delta t^{3}f_{3}+O(\Delta t^{4}) \right]\right)=0 \, .
\end{equation*}
Since the exact Fokker--Planck operator preserves the invariant canonical distribution, i.e., $\mathcal{L}^{\dag}\tilde{\rho}_{\beta}=0$, we obtain
\begin{equation}\label{eq:error_analysis_PDE}
  \mathcal{L}^{\dag}(\tilde{\rho}_{\beta}f_{1}) = - \mathcal{L}^{\dag}_{1}\tilde{\rho}_{\beta}
\end{equation}
by equating first order terms in $\Delta t$.

For any particular integration scheme it is possible to find the perturbation operator $\mathcal{L}^{\dag}_{1}$ by using the BCH expansion.  Then we can calculate its action on $\tilde{\rho}_{\beta}$.  The last step, namely obtaining the leading correction function $f_{1}$, requires the solution of the above PDE (see examples in Langevin dynamics \cite{Leimkuhler2013}).  In general, solving for $f_1$ in closed form is difficult, and it does not get simpler as we consider, as here, more complicated formulations than Langevin dynamics and more complicated splittings.

According to the BCH expansion, for (noncommutative) linear operators $X$ and $Y$, we have
\begin{equation*}
  \exp(\Delta t X)\exp(\Delta t Y) = \exp(\Delta t Z_{1}) \, ,
\end{equation*}
where
\begin{equation}
  Z_{1} = X+Y + \frac{\Delta t}{2}[X,Y] + \frac{\Delta t^{2}}{12}([X,[X,Y]]-[Y,[X,Y]]) + O(\Delta t^{3}) \, ,
\end{equation}
and subsequently
\begin{equation*}
  \exp\left(\frac{\Delta t}{2}X\right)\exp(\Delta t Y)\exp\left(\frac{\Delta t}{2}X\right) = \exp(\Delta t Z_{2}) \, ,
\end{equation*}
where
\begin{equation}
  Z_{2} = X+Y + \frac{\Delta t^{2}}{12}\left([Y,[Y,X]]-\frac{1}{2}[X,[X,Y]]\right) + O(\Delta t^{4}) \, .
\end{equation}
The notation $[X,Y]=XY-YX$ denotes the commutator of operators $X$ and $Y$.

These equations demonstrate that for nonsymmetric splitting methods, there typically exists a nonzero term $\mathcal{L}^{\dag}_{1} \propto [X,Y] \neq 0$, while the condition $\mathcal{L}^{\dag}_{1} = 0$, implying $f_{1}=0$, is automatically satisfied for symmetric splitting methods; thus, for observables $\phi(\q,\p,\xi)$, assuming the asymptotic expansion holds,  the computed average would be of order two
\begin{equation*}
  \langle\phi\rangle_{\Delta t} = \langle\phi\rangle + \Delta t\langle\phi f_{1}\rangle + \Delta t^{2}\langle\phi f_{2}\rangle + \dots = \langle\phi\rangle + O(\Delta t^{2}) \, ,
\end{equation*}
where $\langle \cdot \rangle$ denotes the average with respect to the target invariant distribution. Therefore, the SGNHT-S method~\eqref{eq:Splitting_SGNHT-S} would have second order convergence for all the observables.

We can work out the leading operator $\mathcal{L}^{\dag}_{1}$ associated with the nonsymmetric \mbox{SGNHT-N}/PAD method~\eqref{eq:Splitting_SGNHT-N} of Ding et al.~\cite{Ding2014},
\begin{equation}
  \mathcal{L}^{\dag}_{1,\mathrm{PAD}} = \frac{1}{2}\left( \left[\mathcal{L}^{\dag}_\mathrm{D}, \mathcal{L}^{\dag}_\mathrm{A}\right] + \left[\mathcal{L}^{\dag}_\mathrm{D}, \mathcal{L}^{\dag}_\mathrm{P}\right] + \left[\mathcal{L}^{\dag}_\mathrm{A}, \mathcal{L}^{\dag}_\mathrm{P}\right] \right) \, .
\end{equation}
It is clear that the leading term $f_{1,\mathrm{PAD}}$ in the perturbed distribution~\eqref{eq:distribution_perturbation}  is in general nonzero. Therefore the nonsymmetric SGNHT-N/PAD method would be expected to exhibit first order convergence to the invariant measure. It should be noted that if certain conditions are satisfied, higher order convergence to the invariant measure would be possible as demonstrated by Abdulle et al.~\cite{Abdulle2014a,Abdulle2014}. However, it can be easily demonstrated that it is not the case here for the SGNHT-N/PAD method. In the presence of a noisy gradient, the Ad-Langevin/SGNHT methods, despite the  stepsize dependency~\eqref{eq:Ad-L-2}, would similarly (and generally) be expected to be first order with respect to the invariant distribution.

\subsection{Superconvergence Property}

Recently, it has been demonstrated in the setting of Langevin dynamics that a particular symmetric splitting method (``BAOAB''), which requires only one force calculation per step, is fourth order for configurational quantities in the ergodic limit and in the limit of large friction~\cite{Leimkuhler2013,Leimkuhler2013c}.

In what follows we demonstrate that the newly proposed SGNHT-S/BADODAB method~\eqref{eq:Splitting_SGNHT-S} effectively inherits the superconvergence property of BAOAB in the setting of Ad-Langevin/\\SGNHT system~\eqref{eq:Ad-L-2} with a clean gradient, in case where the parameters $\sigma_{\mathrm{A}}$ and $\mu$ are both taken to infinity in a suitable way.   For simplicity, we consider here a one-dimensional model $H=p^2/2+U(q)$, but the analysis could easily be extended to higher dimensions.

Following the standard procedure described in Section~\ref{subsec:Order_of_Convergence}, we obtain the following PDE associated with the BADODAB method:
\begin{equation}\label{eq:error_analysis_PDE_Ad-L}
  \mathcal{L}^{\dag}(\tilde{\rho}_{\beta}f_{2}) = - \mathcal{L}^{\dag}_{2}\tilde{\rho}_{\beta} \, ,
\end{equation}
where $\mathcal{L}^{\dag}$ is the exact Fokker--Planck operator
\begin{equation}
  \mathcal{L}^{\dag} = - p\partial_{q} + U'(q)\partial_{p} + \xi\partial_{p} (p\cdot) + \frac{\hat{\gamma}}{\beta}\partial_{pp} - \frac{1}{\mu}(p^{2}-\beta^{-1})\partial_{\xi}
\end{equation}
with invariant measure
\begin{equation}
  \tilde{\rho}_{\beta}(q,p,\xi) = \frac{1}{Z}\exp\left({-\beta H(q,p)}\right)\exp\left( -\frac{\beta\mu}{2} (\xi - \hat{\gamma})^{2} \right) \, ,
\end{equation}
where $\hat{\gamma} = \langle \xi \rangle =\beta \sigma^{2}_{\mathrm{A}}/2$ and $\mathcal{L}^{\dag}_{2}$ can be calculated by using the BCH expansion
\begin{align*}
  \mathcal{L}^{\dag}_{2} = \, & \frac{1}{12}\left( \left[\mathcal{L}^{\dag}_\mathrm{O}, \left[\mathcal{L}^{\dag}_\mathrm{O}, \mathcal{L}^{\dag}_\mathrm{D}\right] \right] + \left[\mathcal{L}^{\dag}_\mathrm{D} + \mathcal{L}^{\dag}_\mathrm{O}, \left[\mathcal{L}^{\dag}_\mathrm{D} + \mathcal{L}^{\dag}_\mathrm{O}, \mathcal{L}^{\dag}_\mathrm{A}\right] \right] + \left[\mathcal{L}^{\dag}_\mathrm{A} + \mathcal{L}^{\dag}_\mathrm{D} + \mathcal{L}^{\dag}_\mathrm{O}, \left[\mathcal{L}^{\dag}_\mathrm{A} + \mathcal{L}^{\dag}_\mathrm{D} + \mathcal{L}^{\dag}_\mathrm{O}, \mathcal{L}^{\dag}_\mathrm{B}\right] \right] \right) \\
  & - \frac{1}{24}\left( \left[\mathcal{L}^{\dag}_\mathrm{D}, \left[\mathcal{L}^{\dag}_\mathrm{D}, \mathcal{L}^{\dag}_\mathrm{O}\right] \right] + \left[\mathcal{L}^{\dag}_\mathrm{A}, \left[ \mathcal{L}^{\dag}_\mathrm{A}, \mathcal{L}^{\dag}_\mathrm{D} + \mathcal{L}^{\dag}_\mathrm{O} \right] \right] + \left[\mathcal{L}^{\dag}_\mathrm{B}, \left[\mathcal{L}^{\dag}_\mathrm{B}, \mathcal{L}^{\dag}_\mathrm{A} + \mathcal{L}^{\dag}_\mathrm{D} + \mathcal{L}^{\dag}_\mathrm{O}\right] \right] \right) \, ,
\end{align*}
whose action on the extended invariant measure reads as
\begin{align*}
  \mathcal{L}^{\dag}_{2}\tilde{\rho}_{\beta}  = \, & \frac{1}{12}\left[ - \beta p^{3}U'''(q) + 4\beta p^{2}\xi^{3} + 3\beta\xi p^{2}U''(q) + 3\beta p U'(q)U''(q) + \frac{6\xi p^{2}}{\mu}\left( 1 - \beta p^{2} \right) \right]\tilde{\rho}_{\beta} \\
  & + \frac{\hat{\gamma}}{12}\left[ 3 U''(q) + 4\xi^{2} - 16\beta p^{2}\xi^{2} - 6\beta U''(q)p^{2} + \frac{6}{\mu}\left( 2\beta p^{4} - 5p^{2}+\beta^{-1} \right) \right]\tilde{\rho}_{\beta} \\
  & + \hat{\gamma}^{2}\xi \left(2\beta p^{2}-1\right)\tilde{\rho}_{\beta} + \hat{\gamma}^{3}\left (\frac{2}{3}-\beta p^{2} \right)\tilde{\rho}_{\beta} \, .
\end{align*}

The equation is very complicated, and we have no direct means of solving it.  However, the additional variable $\xi$ has mean $\hat{\gamma}$.  If we suppose that $\mu$ is large, then the variance of $\xi$ will be small.  In this case we can consider the approximation obtained by replacing
functions of $\xi$ in the PDE~\eqref{eq:error_analysis_PDE_Ad-L} by their corresponding averages
\begin{equation}
  \langle\xi\rangle = \hat{\gamma} \, , \qquad \langle\xi^{2}\rangle = \frac{1}{\beta\mu} + \hat{\gamma}^{2} \, , \qquad \langle\xi^{3}\rangle = \frac{3\hat{\gamma}}{\beta\mu} + \hat{\gamma}^{3} \, .
\end{equation}

We use this as part of an averaging of the stationary Fokker--Planck equation with respect to the auxiliary variable.  That is, we project the Fokker--Planck equation and its solution by integrating with respect to the Gaussian distribution of $\xi$ in the ergodic limit.   We can think of this is as defining a sort of ``subspace projection''; it is related to the Galerkin method that is widely used in solving high-dimensional linear systems and PDEs, including Fokker--Planck equations~\cite{Risken1989,Chakravorty2006}.   In this case, we apply the projection operator~\cite{Givon2004}
\begin{equation}
  \mathcal{P} \nu(q,p,\xi) := \frac{\int_{\Omega_\xi} \tilde{\rho}_{\beta}(q,p,\xi) \nu(q,p,\xi) \, {\rm d} \xi}{\int_{\Omega_\xi} \tilde{\rho}_{\beta}(q,p,\xi) \, {\rm d} \xi} \, ,
\end{equation}
where $\nu$ is an arbitrary function, to the PDE~\eqref{eq:error_analysis_PDE_Ad-L}.  Effectively, this results in the reduced equation
\begin{equation}\label{eq:error_analysis_PDE_Ad-Lreduced}
  \check{\mathcal{L}}^{\dag}({\rho}_{\beta}\hat{f}_{2}) = - \rho_{\beta} \mathcal{P} \frac{\mathcal{L}^{\dag}_{2}\tilde{\rho}_{\beta}}{\tilde{\rho}_{\beta}} \, ,
\end{equation}
where the operator $\check{\mathcal{L}}^{\dag}$ is just the operator $\mathcal{L}^{\dag}$ reduced by the action of the projection, and which acts on functions of $q$ and $p$; this is nothing other than the corresponding adjoint generator of Langevin dynamics. Likewise, $\hat{f}_2$ is now a function of $q$ and $p$ only.
The right-hand side simplifies to
\begin{equation*}
 \rho_{\beta} \mathcal{P} \frac{\mathcal{L}^{\dag}_{2}\tilde{\rho}_{\beta}}{\tilde{\rho}_{\beta}} = \left( \frac{\beta}{12}\left[ 3 p U'(q)U''(q) - p^{3}U'''(q) \right] + \frac{\hat{\gamma}}{12}\left[ 3 U''(q) - 3 \beta p^{2}U''(q) + \frac{1}{\mu}\left( 6\beta p^{4} - 28 p^{2} + 10\beta^{-1} \right) \right] \right) \rho_{\beta} \, ,
\end{equation*}
where $\rho_{\beta}$ is the Gibbs (canonical) density (${\rm exp}(-\beta H(q,p))$).

We consider the high friction limit ($\hat{\gamma} \rightarrow \infty$) and expand $\hat{f}_{2}$ in a series involving the reciprocal friction $\varepsilon = 1/\hat{\gamma}$,
\begin{equation}
  \hat{f}_{2}(q,p) = \hat{f}_{2,0}(q,p) + \varepsilon \hat{f}_{2,1}(q,p) + \varepsilon^{2} \hat{f}_{2,2}(q,p) + \cdots \, ,
\end{equation}
with each function $\hat{f}_{2,i}$ satisfying $\langle \hat{f}_{2,i} \rangle=0$. Dividing~\eqref{eq:error_analysis_PDE_Ad-L} by the friction coefficient $\hat{\gamma}$, we obtain
\begin{equation}
  \left( \bar{\mathcal{L}}^{\dag}_\mathrm{O} + \varepsilon \mathcal{L}^{\dag}_\mathrm{H} \right) \left( \hat{f}_{2,0} + \varepsilon \hat{f}_{2,1} + O(\varepsilon^{2}) \right)\rho_{\beta} = - \varepsilon \rho_{\beta} \mathcal{P} \frac{\mathcal{L}^{\dag}_{2}\tilde{\rho}_{\beta}}{\tilde{\rho}_{\beta}} \, ,
\end{equation}
where
\begin{equation}
  \bar{\mathcal{L}}^{\dag}_\mathrm{O} = \partial_{p} (p\cdot) + \beta^{-1}\partial_{pp} \, , \qquad \mathcal{L}^{\dag}_\mathrm{H} = - p\partial_{q} + U'(q)\partial_{p} \, .
\end{equation}
We take the high thermal mass limit ($\mu \rightarrow \infty$) in such a way that $\varepsilon = 1/\mu=1/\hat{\gamma}$. The use of this limit yields the following terms of the expansion of the right-hand side in powers of $\varepsilon$. Defining
\[
- \varepsilon \rho_{\beta} \mathcal{P} \frac{\mathcal{L}^{\dag}_{2}\tilde{\rho}_{\beta}}{\tilde{\rho}_{\beta}} \equiv g =  \left( g_{0} + \varepsilon g_{1}  \right ) \rho_{\beta}\, ,
\]
we have
\begin{align}
  g_{0} & = - \frac{1}{4}\left[ U''(q) - \beta p^{2}U''(q) \right] \, ,\\
  g_{1} & = - \frac{1}{12}\left[ 3\beta p U'(q)U''(q) - \beta p^{3}U'''(q) + 6\beta p^{4} - 28 p^{2} + 10\beta^{-1} \right] \, .
\end{align}
Furthermore, by equating powers of the reciprocal friction $\varepsilon$, we can solve a sequence of equations
\begin{align*}
  & \, \bar{\mathcal{L}}^{\dag}_\mathrm{O}(\rho_{\beta}\hat{f}_{2,0}) = g_{0}\rho_{\beta} \, , \\
  \mathcal{L}^{\dag}_\mathrm{H}(\rho_{\beta}\hat{f}_{2,0}) \, + & \, \bar{\mathcal{L}}^{\dag}_\mathrm{O}(\rho_{\beta}\hat{f}_{2,1}) = g_{1}\rho_{\beta} \, , \\
  \mathcal{L}^{\dag}_\mathrm{H}(\rho_{\beta}\hat{f}_{2,1}) \, + &  \, \bar{\mathcal{L}}^{\dag}_\mathrm{O}(\rho_{\beta}\hat{f}_{2,2}) = 0 \, , \\
  & \vdots
\end{align*}
to obtain the leading term $\hat{f}_{2,0}$, i.e.,
\begin{equation}
  \hat{f}_{2,0} \equiv \hat{f}^{\mathrm{BADODAB}}_{2,0} = \frac{1}{8}\left( U''(q) - \beta p^{2}U''(q) \right) \, .
\end{equation}

Moreover, it can be easily shown that the marginal average of $\hat{f}^{\mathrm{BADODAB}}_{2,0}$ with respect to momentum is zero, i.e.,
\begin{equation}
  \int_{\Omega_p} \hat{f}^{\mathrm{BADODAB}}_{2,0}(q,p) \rho_{\beta} \, \dd \omega_p = 0 \, ,
\end{equation}
which leads to the average of configurational observables $\phi(q)$ with respect to the invariant measure as
\begin{equation*}
  \langle\phi(q)\rangle_{\mathrm{BADODAB}} = \langle\phi(q)\rangle +  \Delta t^{2}\langle\phi(q)\hat{f}^{\mathrm{BADODAB}}_{2,0} \rangle + O(\varepsilon \Delta t^{2} + \Delta t^{4}) \, .
\end{equation*}
Thus, for configurational observables the BADODAB method has fourth order convergence to the invariant measure in the large friction and thermal mass limits (i.e., $\varepsilon \rightarrow 0$),
\begin{equation*}
  \lim_{\varepsilon \rightarrow 0} \langle\phi(q)\rangle_{\mathrm{BADODAB}} = \langle\phi(q)\rangle + O(\Delta t^{4}) \, .
\end{equation*}

It should be emphasized here that only the BADODAB and BAODOAB methods appear to have the superconvergence property among a number of different splitting methods investigated in the Ad-Langevin/SGNHT system~\eqref{eq:Ad-L-2} with a clean gradient. The superconvergence property suggests the use of relatively large $\sigma_{\mathrm{A}}$ and $\mu \propto \sigma^{2}_{\mathrm{A}}$ in the BADODAB (SGNHT-S) method in order to enhance sampling accuracy.  In fact, we expect that larger values of $\mu$ than this bound will not diminish the sampling accuracy, but the effect of large values of $\mu$ is to reduce the responsiveness of the thermostat device.

\section{Numerical Experiments}
\label{sec:Numerical_Experiments}

In this section, we conduct a variety of numerical experiments to compare the performance of the different schemes presented in this article.

\subsection{Molecular Systems}

\begin{figure}[tb]
\centering
\includegraphics[scale=0.5]{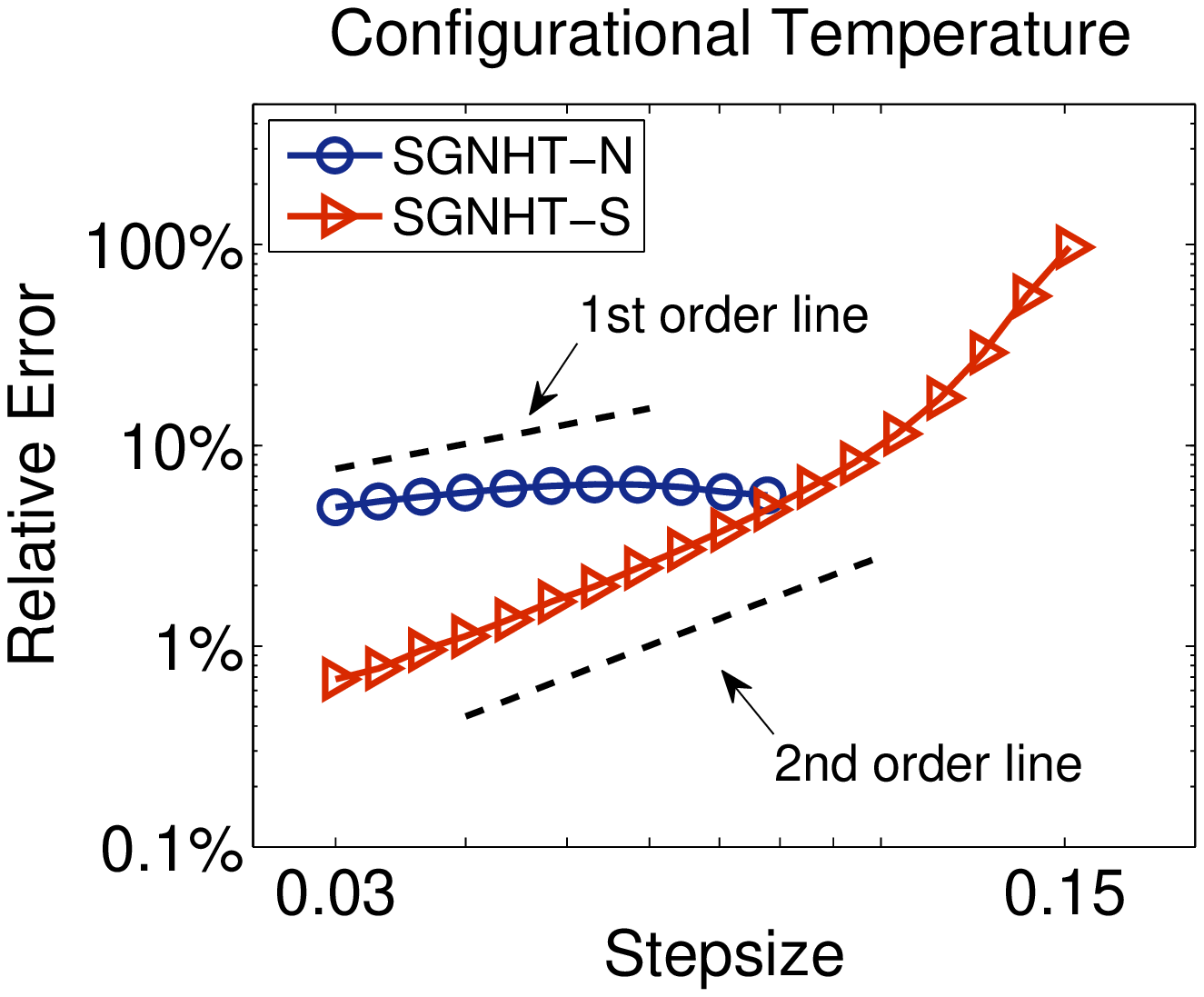}
\includegraphics[scale=0.5]{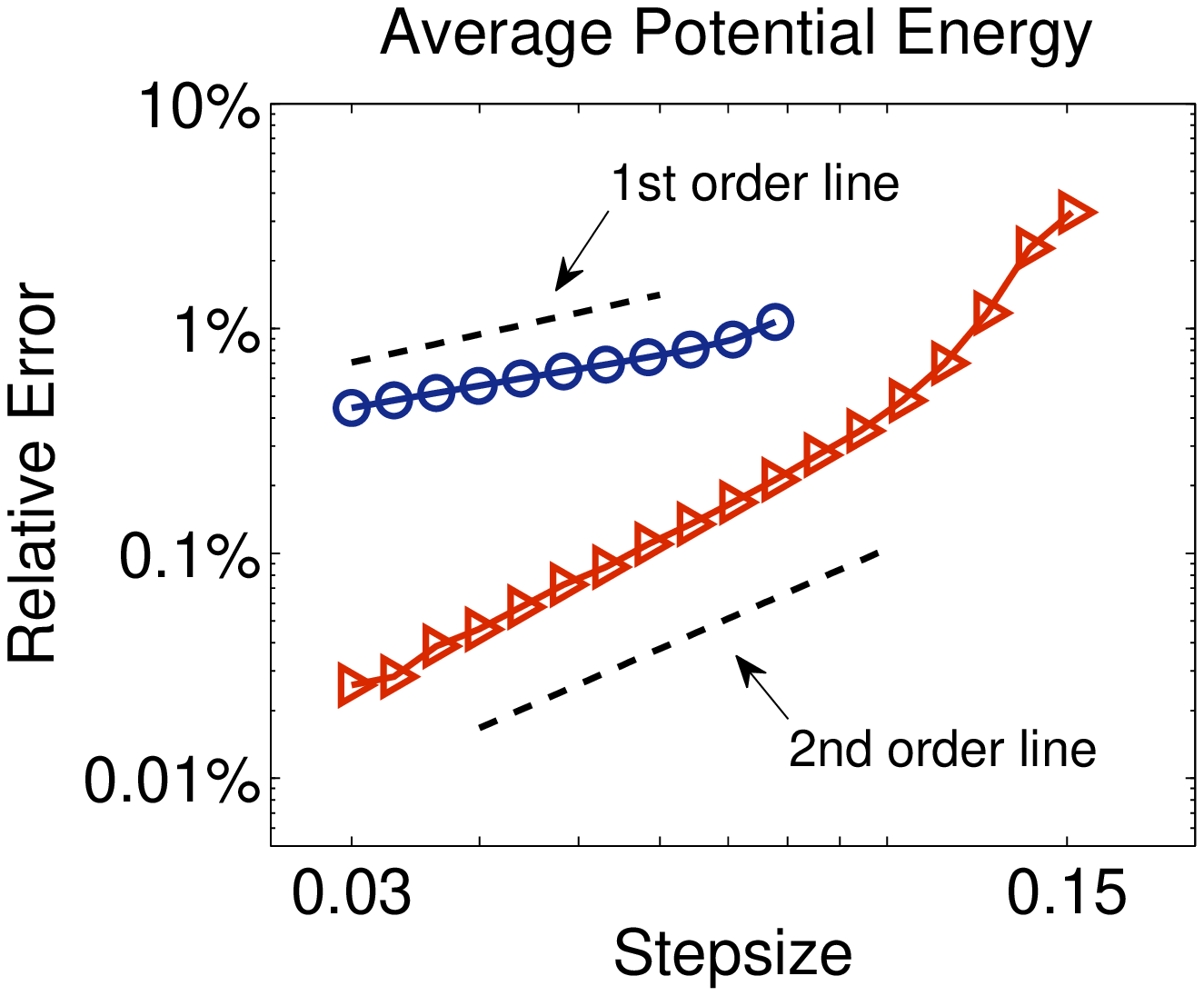}
\caption{\small Log-log plot of the relative error in computed configurational temperature (left) and average potential energy (right) against stepsize by using two Ad-Langevin/SGNHT methods (with a clean gradient). The system ($\sigma_{\mathrm{A}}=3$) was simulated for 5000 reduced time units, but only the last 80\% of the data were collected to calculate the quantity to make sure the system was well equilibrated. Ten different runs were averaged to further reduce the sampling errors. The stepsizes tested began at $\Delta t=0.03$ and were increased incrementally by 10\% until both methods showed significant relative error (SGNHT-N became unstable at around $\Delta t=0.08$). }
\label{fig:Ad-L_Comp_CT_U}
\end{figure}

\begin{figure}[tb]
\centering
\includegraphics[scale=0.5]{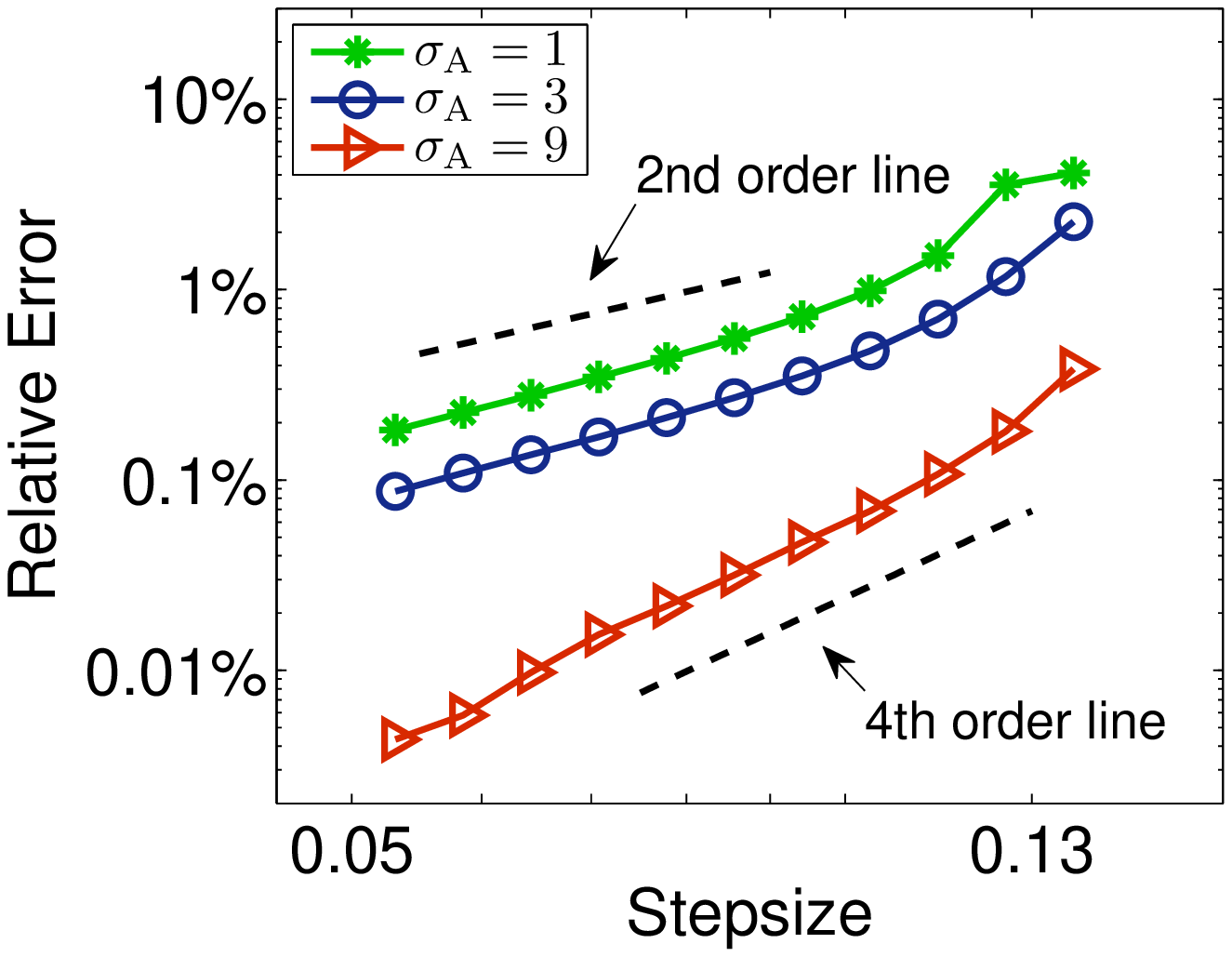}
\includegraphics[scale=0.5]{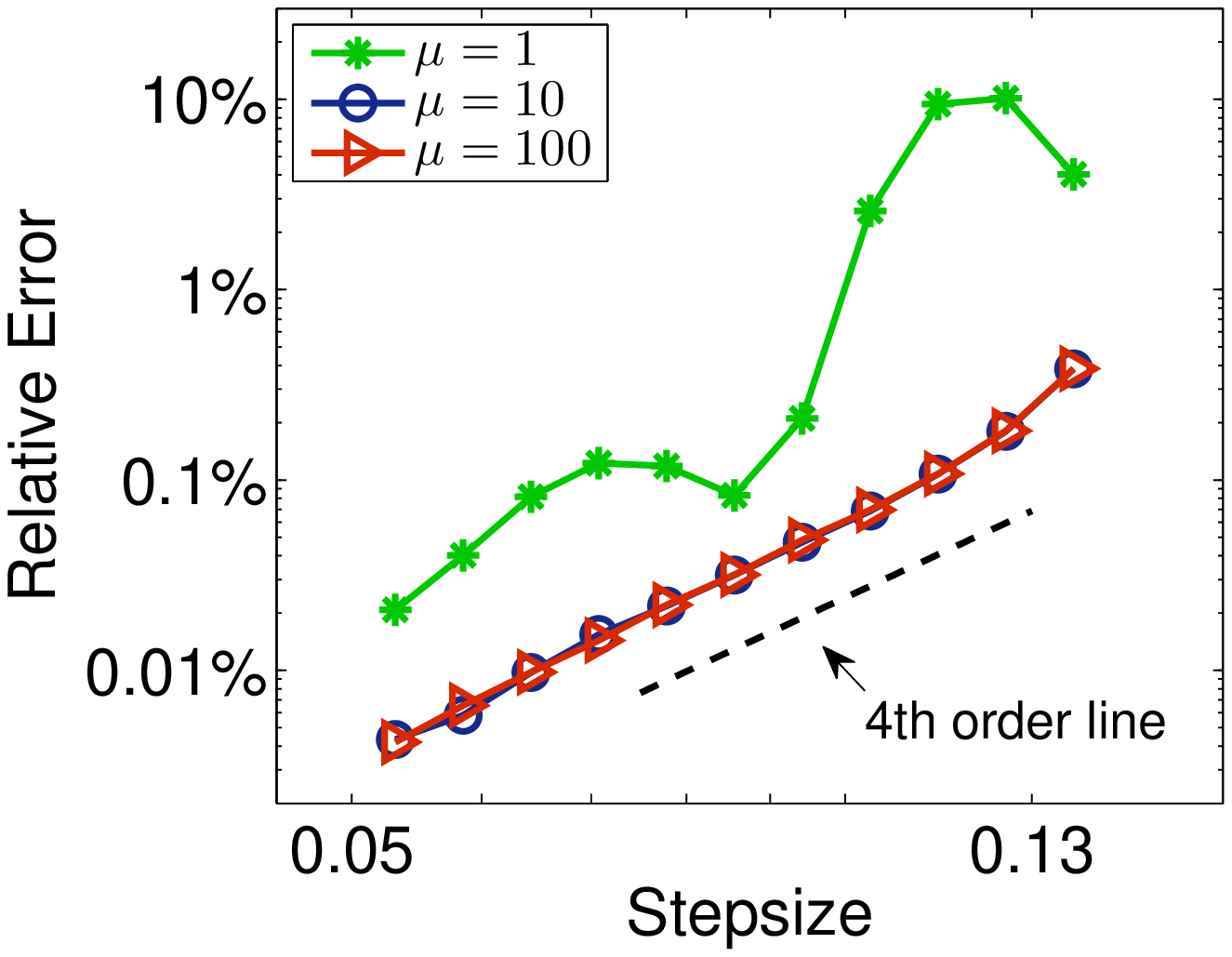}
\caption{\small Log-log plot of the relative error in computed average potential energy against stepsize by using the SGNHT-S/BADODAB method with (left) different values of $\sigma_{\mathrm{A}}$ ($\mu=10$)   and (right) different values of $\mu$ ($\sigma_{\mathrm{A}}=9$). The format of the plots is the same as in Figure~\ref{fig:Ad-L_Comp_CT_U} except 50 different runs were used to reduce the sampling errors in high accuracy regime. }
\label{fig:Ad-L_Comp_sigma_mu}
\end{figure}

Before we compare various methods in machine learning applications (i.e., with a noisy gradient), we first demonstrate the order of convergence of various splitting methods with a clean gradient.

A popular model of an $N$-body system with pair interactions based on a spring with rest length (i.e., pendulum) was used, a standard if simplified model of molecular dynamics. The total potential energy of the system is defined as
\begin{equation}
  U(\q) = \sum^{N-1}_{i=1} \sum^{N}_{j=i+1}\varphi(r_{ij}) \, ,
\end{equation}
where $r_{ij}=\|\q_{i}-\q_{j}\|$ denotes the distance between two particles $i$ and $j$, and $\varphi(r_{ij})$ represents the pair potential energy
\begin{equation}\label{eq:Potential_Pendulum}
  \varphi(r_{ij})=
  \begin{cases}
    \displaystyle\frac{k}{2}\left( r_{ij} - r_{\mathrm{c}} \right)^{2} \, , & r_{ij}<r_{\mathrm{c}} \, ;\\
    \quad \quad \ 0 \, , & r_{ij}\geq r_{\mathrm{c}} \, ,
  \end{cases}
\end{equation}
where $k$ and $r_{\mathrm{c}}$ represent the spring constant and the cutoff radius, respectively.

A system consisting of $N=500$ identical particles (i.e., unit mass) was simulated in a cubic box with periodic boundary conditions~\cite{Allen1989}. The positions of the particles were initialized on a cubic grid with equidistant grid spacing, while the initial momenta were i.i.d.\ random variables with mean zero and variance $k_{\mathrm{B}}T$, which was set to be unity. The thermal mass $\mu$ was chosen to be 10 unless otherwise stated. Particle density $\rho_{\rm d}=4$ was used with spring constant $k=25$ and cutoff radius $r_{\mathrm{c}}=1$.

We first compare the two SGNHT methods on controlling two configurational quantities: configurational temperature and average potential energy. The configurational temperature~\cite{Hirschfelder1960}, which, as the kinetic temperature, should in principle be equal to the target temperature, can be defined as
\begin{equation*}
  k_{\mathrm{B}}T = \frac{\sum_{i}\langle {\| \nabla_{i}U \|}^{2} \rangle}{\sum_{i}\langle \nabla^{2}_{i}U \rangle} \, ,
\end{equation*}
where the angle brackets denote the averages, and $\nabla_{i}U$ and $\nabla^{2}_{i}U$ represent the gradient and Laplacian of the potential energy $U$ with respect to the position of particle $i$, respectively (see more discussions in~\cite{Leimkuhler2015}).

As shown in Figure~\ref{fig:Ad-L_Comp_CT_U}, with the help of the dashed order lines, we can see that SGNHT-N and SGNHT-S show first and second order convergence, respectively, as expected. It is clear that SGNHT-S has not only at least one order of magnitude improvement in accuracy in both observables, but also much greater robustness over the SGNHT-N method, which becomes completely unstable at around $\Delta t=0.08$. The results on the configurational temperature and average potential energy are rather similar; therefore in what follows we present only average potential energy results.

We also investigate the effect of changing the value of $\sigma_{\mathrm{A}}$ in the SGNHT-S/BADODAB scheme proposed in this article. As can be seen from Figure~\ref{fig:Ad-L_Comp_sigma_mu}, the SGNHT-S method displays second order convergence to the invariant measure when $\sigma_{\mathrm{A}}$ is relatively small, while a fourth order convergence is observed in the high friction limit ($\sigma_{\mathrm{A}}=9$), as anticipated from the analysis of the previous section. It should be emphasized here that the superconvergence property was observed only in the BADODAB and BAODOAB methods, which both reduce to the BAOAB method~\cite{Leimkuhler2013,Leimkuhler2013c} in Langevin dynamics.

Figure~\ref{fig:Ad-L_Comp_sigma_mu} also compares the effect of varying the value of the thermal mass $\mu$ when $\sigma_{\mathrm{A}}$ is fixed. It can be seen that the BADODAB method displays a clear fourth order convergence when $\mu$ is relatively large, while when $\mu$ is small, not only is the smooth discretization error dependence on stepsize lost, but significantly larger relative error is also observed. This reinforces the choice of a relatively large value of $\mu$. It is worth pointing out that $\mu=10$ works as well as $\mu=100$; therefore $\mu=10$ is used throughout this article since a relatively smaller $\mu$ corresponds to a tighter interaction between the thermostat and the system, and thus it can fluctuate more rapidly to accommodate changes in the noise and adapt more easily.

\begin{figure}[tb]
\centering
\includegraphics[scale=0.5]{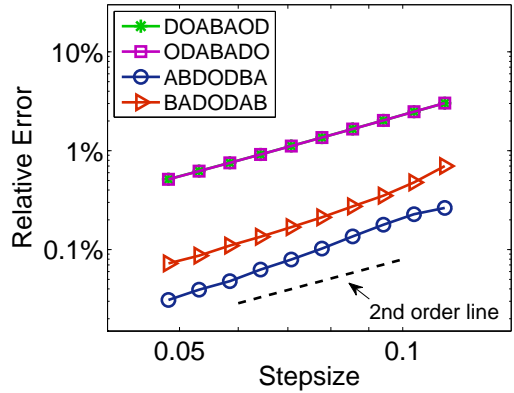}
\caption{\small Log-log plot of the relative error in computed average potential energy against stepsize by using various splitting methods of the Ad-Langevin/SGNHT system ($\sigma_{\mathrm{A}}=3$). The format of the plot is the same as in Figure~\ref{fig:Ad-L_Comp_CT_U}. }
\label{fig:Ad-L_Comp_Methods}
\end{figure}

We also explore in Figure~\ref{fig:Ad-L_Comp_Methods} the performance of various splitting methods of the \mbox{Ad-Langevin}/\\SGNHT system~\eqref{eq:Ad-L-2} with fixed values of $\sigma_{\mathrm{A}}$ and $\mu$. All the methods clearly show second order convergence, with ABDODBA and BADODAB methods achieving one order of magnitude improvement in accuracy compared to the other methods. This again illustrates the importance of optimal design of numerical methods. The ABDODBA method seems to be slightly better that the BADODAB method in the regime of $\sigma_{\mathrm{A}}=3$; however, as demonstrated in Figure~\ref{fig:Ad-L_Comp_sigma_mu}, the BADODAB method achieves a dramatic improvement in accuracy when $\sigma_{\mathrm{A}}$ is relatively large (e.g., $\sigma_{\mathrm{A}}=9$), while other schemes remain the same except for the BAODOAB method.

\subsection{Bayesian Inference}
\label{subsec:Bayesian_Inference}

In this subsection we compare methods in a classical Bayesian inference model in one dimension, i.e., to estimate the mean of a normal distribution with known variance~\cite{Ding2014}. More precisely, given $N$ i.i.d.\ samples from a normal distribution, $x_{i} \sim \mathcal{N}(\check{\mu},\hat{\sigma}^{2})$, where it should be noted that $\check{\mu}$ is the true mean, when we draw samples with known $\hat{\sigma}^{2}$ and a uniform prior distribution ranging from $-N/2$ to $N/2$, we are able to calculate the posterior distribution of the mean in a closed form
\begin{equation}\label{eq:Bayesian_Inf_Post_Mean_1}
  \hat{\mu} \sim \mathcal{N}\left(\hat{x},\frac{\hat{\sigma}^{2}}{N}\right) \, ,
\end{equation}
where $\hat{x}=\sum^{N}_{i=1}x_{i}/N$. In the context of stochastic gradient approximation, we have
\begin{equation}
  \label{eq:Bayesian_Inf_Post_Mean_2}
  \begin{aligned}
    \pi(\hat{\mu}|\mathbf{X}) &\propto \pi(\mathbf{X}|\hat{\mu}) \pi(\hat{\mu}) \approx  \left( \prod_{i=1}^{\tilde{n}} \pi(\mathbf{x}_{r_{i}}|\hat{\mu}) \right)^{\frac{N}{\tilde{n}}} \pi(\hat{\mu}) \\
    &= \left( \frac{1}{\sqrt{2\pi}\hat{\sigma}} \right)^{N} \left[ \prod_{i=1}^{\tilde{n}}\exp \left( - \frac{(x_{i}-\hat{\mu})^{2}}{2\hat{\sigma}^{2}}\right) \right]^{\frac{N}{\tilde{n}}} \frac{1}{N} \\
    &= \left( \frac{1}{\sqrt{2\pi}\hat{\sigma}} \right)^{N} \exp \left( - \frac{N}{\tilde{n}} \sum_{i=1}^{\tilde{n}} \frac{(x_{i}-\hat{\mu})^{2}}{2\hat{\sigma}^{2}}\right) \frac{1}{N} \\
    &\propto \exp \left( - \frac{N}{\tilde{n}} \sum_{i=1}^{\tilde{n}} \frac{(x_{i}-\hat{\mu})^{2}}{2\hat{\sigma}^{2}}\right) \\
    &= \exp \left[ - \frac{1}{2\hat{\sigma}^{2}} \frac{N}{\tilde{n}} \left( \sum_{i=1}^{\tilde{n}} (x_{i}-\bar{x})^{2} + \tilde{n}(\bar{x}-\hat{\mu})^{2} \right) \right] \\
    &\propto \exp \left( - \frac{N}{2\hat{\sigma}^{2}} (\bar{x}-\hat{\mu})^{2} \right) \, ,
  \end{aligned}
\end{equation}
where $\bar{x}=\sum^{\tilde{n}}_{i=1}x_{i}/\tilde{n}$. It clearly recovers the true distribution~\eqref{eq:Bayesian_Inf_Post_Mean_1} when $\tilde{n}=N$. Taking the logarithm and differentiating the posterior distribution obtained at the end of~\eqref{eq:Bayesian_Inf_Post_Mean_2} with respect to $\hat{\mu}$ gives the noisy force
\begin{equation}\label{eq:Bayesian_Inf_Noisy_F}
  \tilde{F}(\hat{\mu}) = \frac{N}{\hat{\sigma}^{2}}\left(\hat{\mu} - \frac{1}{\tilde{n}}\sum^{\tilde{n}}_{i=1}x_{i} \right) \, .
\end{equation}

\begin{figure}[tb]
\centering
\includegraphics[scale=0.5]{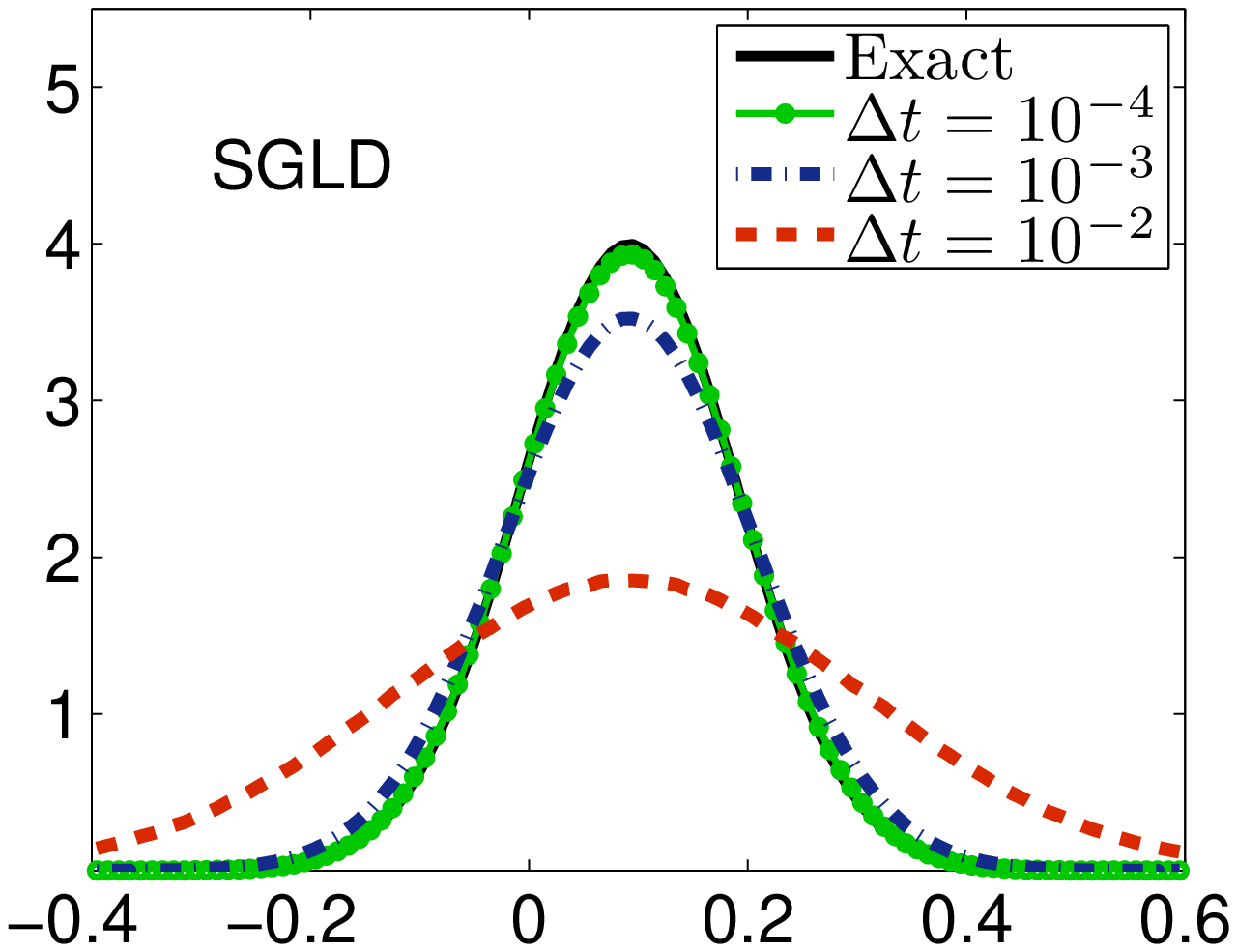}
\includegraphics[scale=0.5]{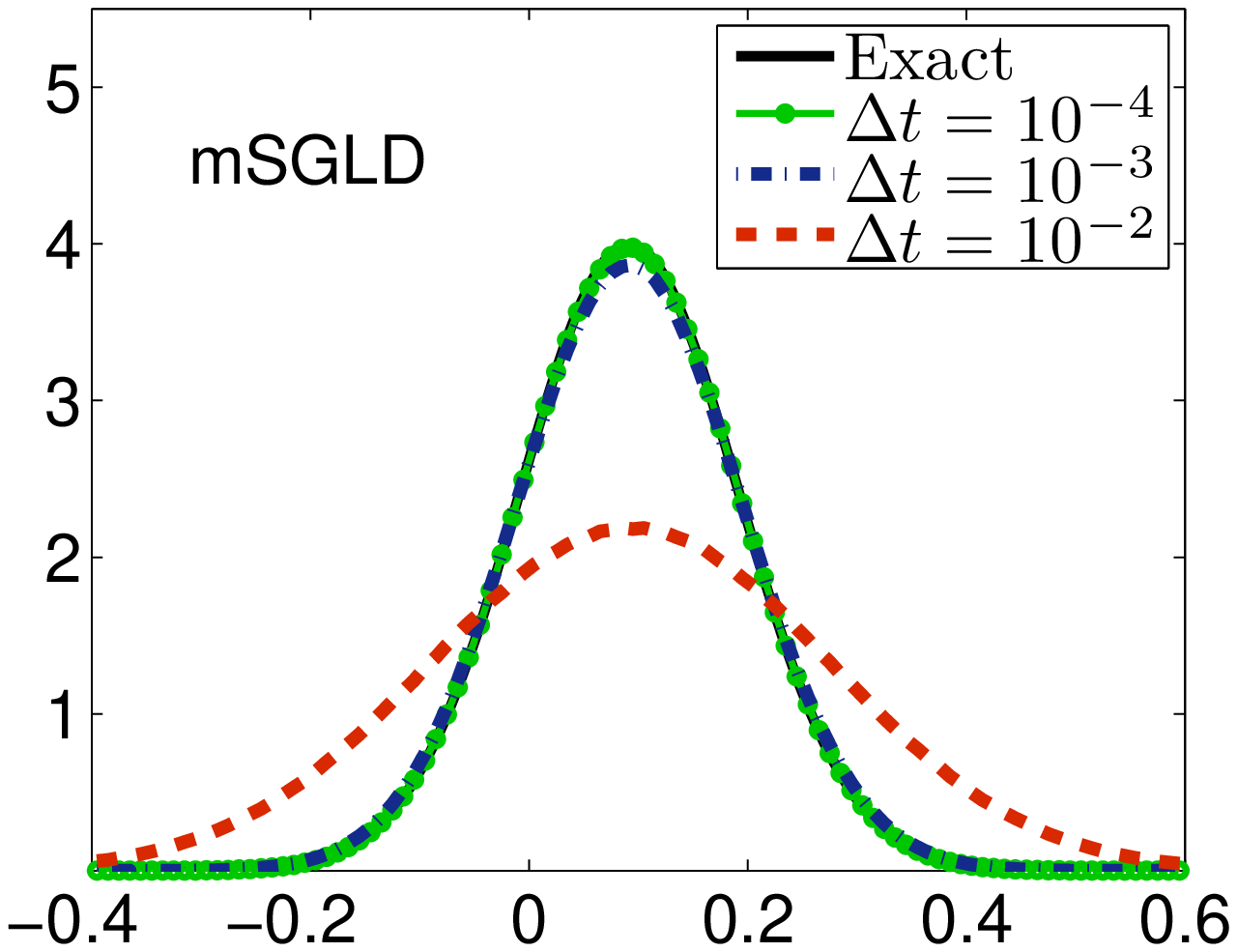}
\includegraphics[scale=0.5]{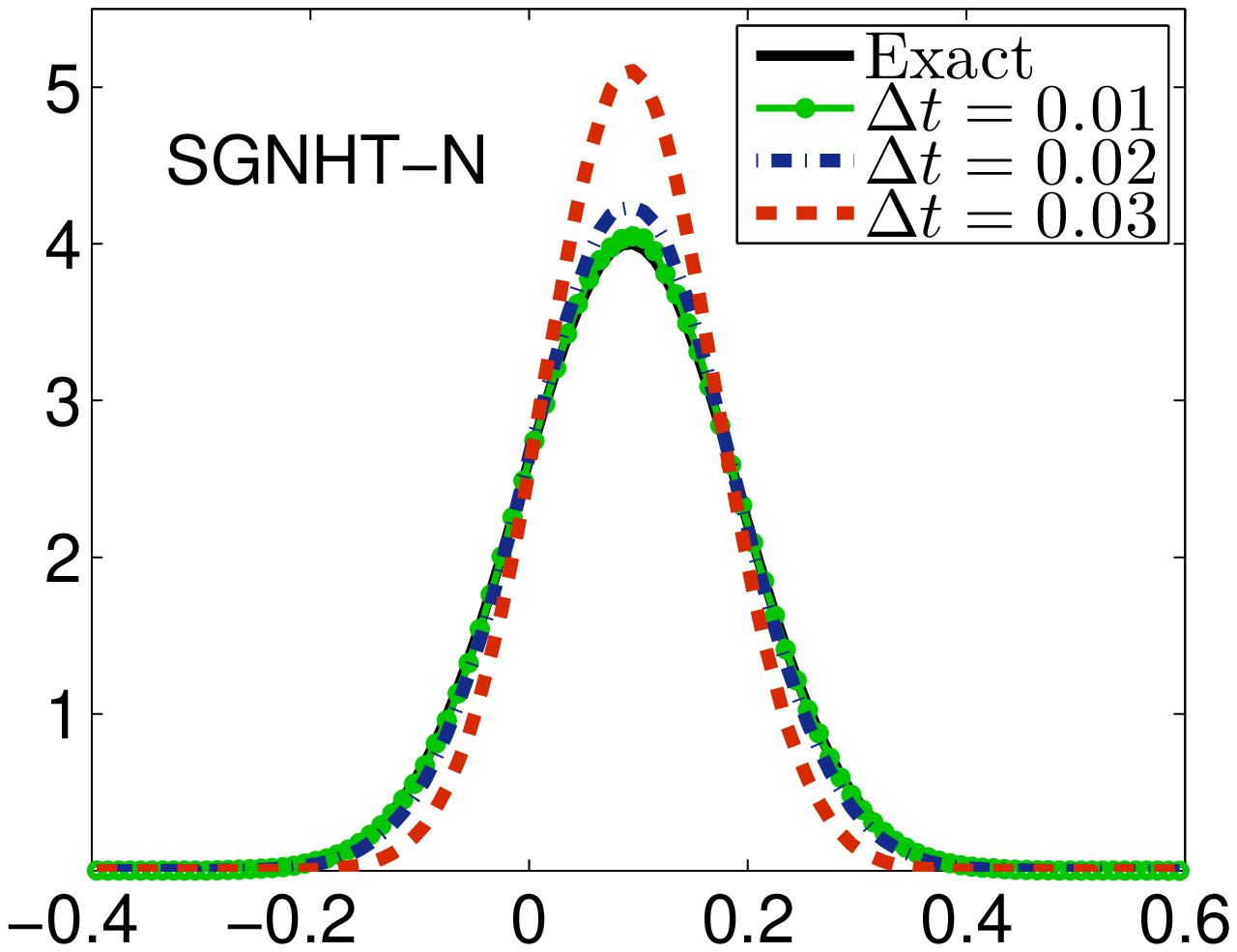}
\includegraphics[scale=0.5]{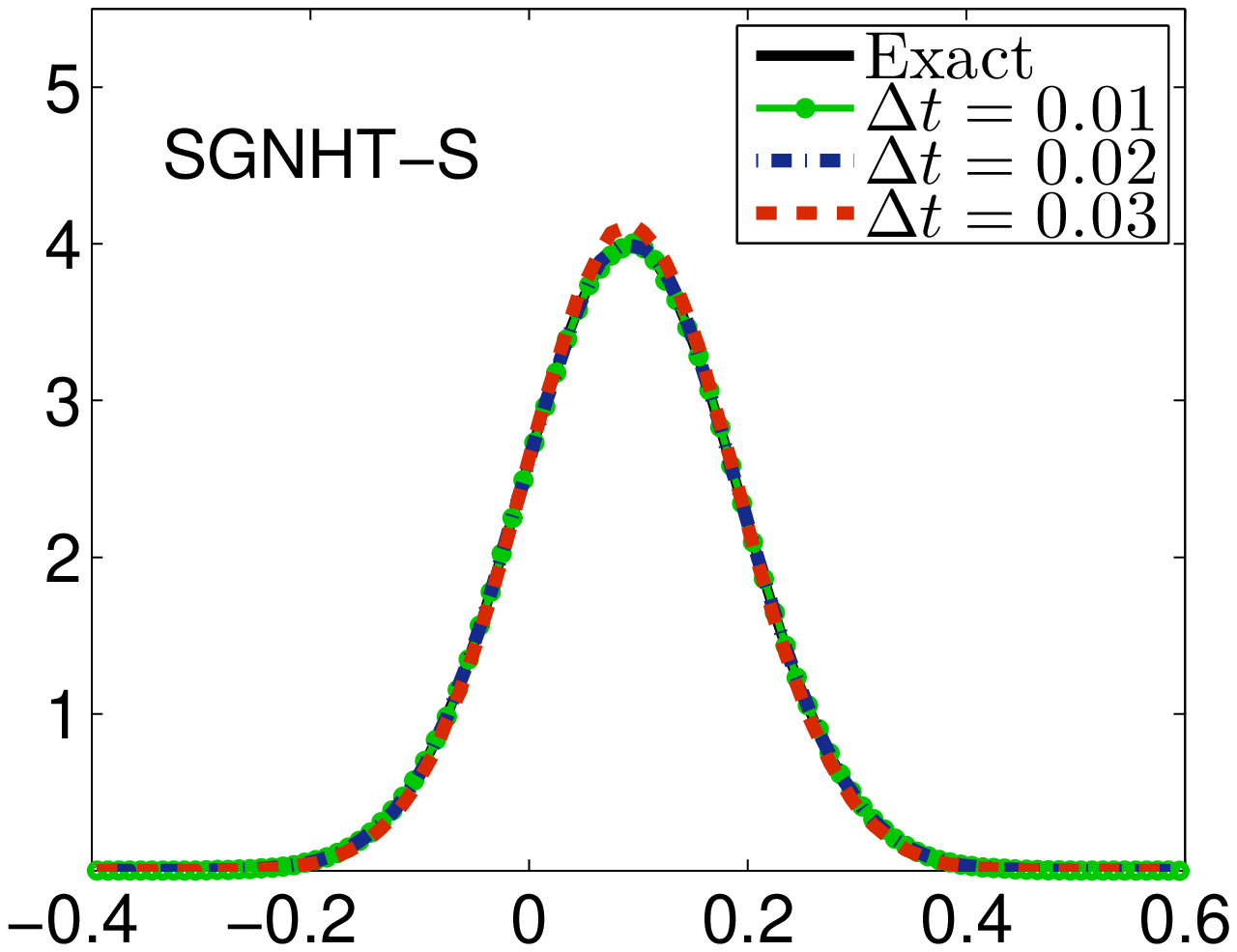}
\caption{\small Comparisons of the distribution in a one-dimensional Bayesian inference problem by using SGLD (top left), mSGLD (top right), SGNHT-N (bottom left), and SGNHT-S (bottom right) with different stepsizes indicated by different colors. The solid black line is the exact solution. Note the difference in the legends between rows. }
\label{fig:Bayes_Inf_Comp_Dist}
\end{figure}

In this simple case, the noise of the stochastic gradient is independent of $\hat{\mu}$ and is a constant given $\tilde{n}$. Moreover, we are able to obtain its mean and variance with respect to the stochastic gradient~\cite{Vollmer2015,Horvitz1952}:
\begin{equation}
  \label{eq:Bayesian_Inf_Mean_Var}
  \begin{aligned}
    \E\tilde{F}(\hat{\mu}) &= F(\hat{\mu}) = \frac{N}{\hat{\sigma}^{2}}\left(\hat{\mu} - \frac{1}{N}\sum^{N}_{i=1}x_{i} \right) \, , \\
    \Var\tilde{F}(\hat{\mu}) &= \frac{1}{\hat{\sigma}^{4}} \frac{N(N-1)}{\tilde{n}} \Var \mathbf{X} \, ,
  \end{aligned}
\end{equation}
where $\Var \mathbf{X}$ is the variance of the dataset. Thus, it is straightforward to verify that the noise is normally distributed according to the central limit theorem.

In our numerical experiments, $\sigma_{\mathrm{A}}$ was chosen as 1 due to the fact that large $\sigma_{\mathrm{A}}$ results in stability issues here. We generated $N=100$ samples from $\mathcal{N}(0,1)$ and randomly selected a subset of size $\tilde{n}=10$ at each timestep to compute the noisy force~\eqref{eq:Bayesian_Inf_Noisy_F}. We plot the distributions of the posterior mean of the dataset obtained by using four different methods with different stepsizes in Figure~\ref{fig:Bayes_Inf_Comp_Dist}. Clearly, two SGNHT methods completely outperformed the SGLD and mSGLD methods. The latter only demonstrate good approximation of the true distribution with order of magnitude smaller stepsize compared to the former. But it should be noted that mSGLD here is slightly better than SGLD in maintaining the true distribution: the distribution of mSGLD with $\Delta t=0.001$ is visibly much closer to the target compared to  that of SGLD with the same stepsize.

Note that stepsizes for SGNHT (second order dynamics) and SGLD (first order dynamics) based methods are not directly comparable---as mentioned in~\cite{Leimkuhler2013} the stepsize of a first order dynamics method like Euler--Maruyama when viewed as the limiting discretization of a Langevin integrator corresponds to $\Delta t^{2}/2$, where $\Delta t$ is the stepsize of the Langevin method.   However, in our experiments we are uninterested in the time-dynamics of the system and care only about the invariant measure.  Therefore the important relationship is the error in thermodynamic averages in comparison with the number of timesteps (work), which quantifies the efficiency of a given method.  The stepsize is just an arbitrary parameter which allows for refinement of the statistical calculation.

\begin{figure}[tb]
\centering
\includegraphics[scale=0.5]{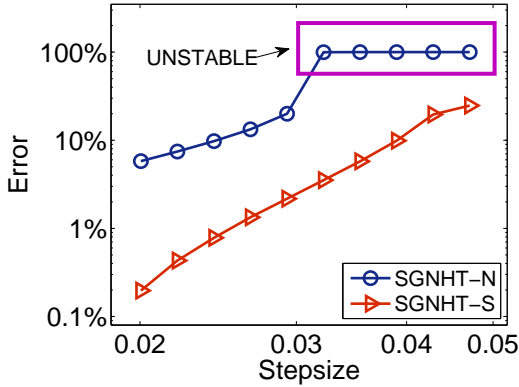}
\caption{\small Log-log plot of the MAE in the distribution of the Bayesian inference model against stepsize. The box indicates that the system was unstable with corresponding stepsizes for the \mbox{SGNHT-N} method. }
\label{fig:Bayes_Inf_Comp_Dist_Error}
\end{figure}

Between the two SGNHT methods, SGNHT-S (the new scheme being proposed here) is obviously superior to SGNHT-N: the latter starts to show significant deviation from the true distribution at $\Delta t=0.02$, while the distribution of the former still looks well matched to the true one at $\Delta t=0.03$. Our observations are confirmed by Figure~\ref{fig:Bayes_Inf_Comp_Dist_Error}, where the mean absolute error (MAE) of the distribution of the two SGNHT methods is plotted. The MAE, which can be thought of as a relative error in distribution, is defined as
\begin{equation}
  \mathrm{MAE} = \frac{1}{\bar{N}} \sum^{\bar{N}}_{i=1} |\omega_{i}-\hat{\omega}_{i}| \, ,
\end{equation}
where $\bar{N}$ denotes the number of intervals, which was chosen as 100. $\omega_{i}$ and $\hat{\omega}_{i}$ represent the observed frequency in bin $i$ and the exact expected frequency, respectively~\cite{Leimkuhler2013}. As can be seen, the stability threshold of SGNHT-N was around $\Delta t=0.03$, beyond which the system became unstable, as highlighted in the figure (in which case the system blew up, resulting in a 100\% MAE). Once again, SGNHT-S not only shows an order of magnitude better accuracy but also has a much greater robustness than SGNHT-N. In particular, for defined accuracy, the SGNHT-S method is able to use double the stepsize compared to SGNHT-N, which means a remarkable 50\% improvement in overall numerical efficiency as defined in~\cite{Leimkuhler2015}.

\subsection{Bayesian Logistic Regression}

Following~\cite{Vollmer2015}, we also investigate the performance of different methods for a more complicated Bayesian logistic regression model. The data $y_i \in \{-1,1\}$ were modelled by
\begin{equation*}
  \pi(y_i | \mathbf{x}_{i}, \boldsymbol{\beta}) = f(y_{i}\boldsymbol{\beta}^{T}\mathbf{x}_{i}) \, ,
\end{equation*}
where $f(z)=1/(1+\exp(-z)) \in [0,1]$ is the logistic function and $\mathbf{x}_{i} \in \mathbb{R}^{d}$ are rows of a fixed dataset. Our goal is to estimate the posterior mean of parameter vector $\boldsymbol{\beta} \in \mathbb{R}^{d}$. For simplicity, a multivariate Gaussian prior $\mathcal{N}(\mathbf{0},\I)$ was used on $\boldsymbol{\beta}$. Therefore, by using Bayes' theorem, we obtain the following posterior distribution:
\begin{equation}
  \pi(\boldsymbol{\beta}) \propto \exp\left(-\frac{1}{2} \|\boldsymbol{\beta}\|^{2}\right) \prod^{N}_{i=1}f(y_{i}\boldsymbol{\beta}^{T}\mathbf{x}_{i}) \, .
\end{equation}
Following the same procedure in the Bayesian inference example (Section~\ref{subsec:Bayesian_Inference}), we can calculate the noisy force and then plug it into different thermostats for sampling.

In our numerical experiments, we considered the $d=3$ case with $N=1000$ data points. We chose the dataset to be
\begin{equation}
\mathbf{X} =
\left(
  \begin{array}{ccc}
    x_{1,1} & x_{1,2} & 1 \\
    x_{2,1} & x_{2,2} & 1 \\
    \vdots & \vdots & \vdots \\
    x_{1000,1} & x_{1000,2} & 1 \\
  \end{array}
\right) \, ,
\end{equation}
where $x_{i,j}$ were sampled from a standard normal distribution $\mathcal{N}(0,1)$ for $i=1,\dots,1000$ and $j=1,2$. A subset of size $\tilde{n}=100$ was randomly chosen at each timestep to compute the noisy force.

\begin{figure}[tb]
\centering
\includegraphics[scale=0.5]{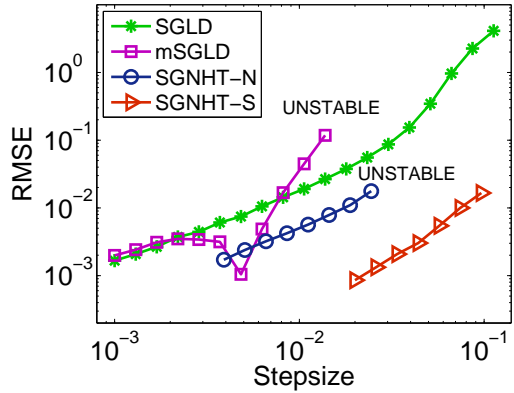}
\caption{\small Comparisons of the RMSE of the posterior mean in the Bayesian logistic regression model by using various methods against stepsize. The system was simulated for 1000 reduced time units with 100,000 different runs. The stepsizes tested began at $\Delta t=0.001$ and were increased incrementally by 30\% until all methods either displayed significant error or became unstable (mSGLD and SGNHT-N). }
\label{fig:Bayes_Logit_Comp_RMSE}
\end{figure}

The performance of estimating the posterior mean value of parameter vector $\boldsymbol{\beta}$ by various methods ($\sigma_\mathrm{A}=6$) was tested and plotted in Figure~\ref{fig:Bayes_Logit_Comp_RMSE}. Again, SGLD and mSGLD, displaying considerably larger root mean square error (RMSE) with a fixed stepsize, were outperformed by the two SGNHT methods.
In this case, the SGLD and mSGLD methods demonstrated similar control in numerical accuracy, but the latter displayed much worse stability than that of the former and became unstable just above $\Delta t=0.01$. As reported in the original paper~\cite{Vollmer2015}, the performance of the mSGLD method depends strongly on the size of the subset---for a larger subset, which requires higher computational cost, the bias of mSGLD can be smaller than that of SGLD.

Of the two SGNHT methods, the SGNHT-S method again shows not only at least an order of magnitude improvement on accuracy but also much better robustness than the other: SGNHT-N became unstable just above $\Delta t=0.02$. Remarkably, the SGNHT-S method at $\Delta t=0.1$ still achieves better accuracy than the SGLD method at $\Delta t=0.01$. In other words, the method we propose here gives more than a 90\% improvement in overall numerical efficiency compared to one of the most popular methods in the literature. For fixed  accuracy, the SGNHT-S method can use almost four times the stepsize of the SGNHT-N method (i.e., an improvement of about 75\% in overall numerical efficiency).

\section{Conclusions}
\label{sec:Conclusions}

We have reviewed a variety of methods in stochastic gradient systems with applications in machine learning. We have provided a theoretical discussion on the foundation (underlying dynamics) of those stochastic gradient systems, which has been lacking in the literature. We have also proposed a new symmetric splitting (at least second order) method in SGNHT (SGNHT-S/BADODAB), which substantially improves the accuracy and robustness compared to a nonsymmetric splitting (first order) method (SGNHT-N) proposed recently in the literature. Furthermore, we have demonstrated that under certain conditions the \mbox{SGNHT-S}/BADODAB method can inherit the superconvergence property recently discovered in integrators for Langevin dynamics, i.e., fourth order convergence  to the invariant measure for configurational averages.

By conducting various numerical experiments, we have demonstrated that the two SGNHT methods outperform the popular SGLD method and its variant mSGLD. In particular, the SGNHT-S method can use up to ten times the stepsize of SGLD, which implies a remarkable more than 90\% improvement in overall numerical efficiency. Between the two SGNHT methods, the SGNHT-S method can use almost four times the stepsize of SGNHT-N for defined accuracy (i.e., about a 75\% improvement in overall numerical efficiency).

It should be noted that in certain cases, it may be desirable to employ a Metropolis--Hastings procedure in order to remove the discretization bias~\cite{Roberts1996}. However, we emphasize that the correction is not without computational cost, particularly as the dimension is increased~\cite{Kennedy1991,Roberts1997,Roberts1998,Beskos2013}, and the results of~\cite{Leimkuhler2013,Leimkuhler2013a,Leimkuhler2013c} and of the current article demonstrate that high accuracy with respect to the invariant distribution is often achievable using traditional numerical integration techniques, thus in many cases entirely eliminating the necessity of Metropolis--Hastings corrections (see more discussions in~\cite{Leimkuhler2013c}). Moreover, we mention that the methods of this article can in principle be combined with Metropolis--Hastings algorithms if it is necessary to completely eliminate the discretization bias.

\section*{Acknowledgements}

The authors thank Ben Goddard, Charles Matthews, Tony Shardlow, Zhanxing Zhu, and Konstantinos Zygalakis for stimulating discussions and valuable suggestions. The authors further thank the anonymous referees for their comments, which substantially contributed to the presentation of our results. XS gratefully acknowledges the financial support from the University of Edinburgh and China Scholarship Council.

\bibliographystyle{is-abbrv}

\bibliography{refs}

\end{document}